\documentclass[draftcls,onecolumn,12pt]{IEEEtran}
\usepackage{amssymb,times}
\usepackage{lipsum}%
\usepackage{cite, epsfig}
\usepackage{verbatim, amsfonts,enumerate}     
\usepackage{bm,array}
\usepackage{graphicx}
\usepackage[cmex10]{amsmath}
\usepackage{epstopdf}
\usepackage[font=small]{caption}
\usepackage{subcaption}
\usepackage{mathtools}					
\usepackage{amsthm}
\usepackage{algorithm,algorithmic}
\usepackage{framed}
\usepackage{color,url}
\usepackage[normalem]{ulem}
\theoremstyle{plain}
\newtheorem{thm}{Theorem}
\newtheorem{lem}{Lemma}

\def\h{\mathbf{h}}
\def\v{\mathbf{v}}
\def\u{\mathbf{u}}
\def\x{\mathbf{x}}
\def\n{\mathbf{n}}
\def\e{\mathbf{e}}
\def\y{\mathbf{y}}
\def\Rn{\mathbb{R}}
\def\E{\mathbf{E}}
\def\X{\mathbf{X}}
\def \Y{\mathbf{Y}}
\def \J{\mathbf{J}}
\def \M {\mathbf{M}}
\def \U {\mathbf{U}}
\def \V {\mathbf{V}}
\def \D {\mathbf{D}}
\def\z{\mathbf{z}}
\def\b{\mathbf{b}}
\def \F {\mathcal{F}}
\def \Fk {\F_{k}}
\def \lb	{{\boldsymbol{\lambda}}}
\def \nb	{{\boldsymbol{\nu}}}

\def \reg {{\mathbf{Reg}}}

\def \cX {{\mathcal{X}}}
\def \Xk {{\cX^\star_k}}
\def \pa {{\mathsf{P_1}}}
\def \pb {{\mathsf{P_2}}}

\newcommand{\colb}[1]{\textcolor{black}{#1}}
\newcommand{\col}[1]{\textcolor{black}{#1}}
\providecommand{\abs}[1]{\left|#1\right|}
\providecommand{\norm}[1]{\left\|#1\right\|}
\providecommand{\ip}[1]{\langle#1\rangle}

\providecommand{\Ex}[1]{\mathbb{E}\left[#1\right]}
\providecommand{\Exc}[1]{\mathbb{E}_{k}\left[#1\right]}

\providecommand{\dist}[1]{\text{dist}(#1)}

\theoremstyle{plain}
\theoremstyle{remark}
\theoremstyle{remark}
\newtheorem{cor}{\bf Corollary} 
\newtheorem{rem}{\bf Remark} 

\graphicspath{{img/}}
\begin{document}


\title{\vspace{-0mm}\colb{Tracking Moving Agents via Inexact Online Gradient Descent Algorithm} \author{Amrit Singh Bedi, \emph{Student Member, IEEE }, Paban Sarma, and Ketan Rajawat, \emph{Member, IEEE }\vspace{0mm}}\thanks{The authors are with the department of Electrical Engineering, Indian Institute of Technology Kanpur, Kanpur, Uttar Pradesh, India -208016.}}

	\maketitle
\begin{abstract}
Multi-agent systems are being increasingly deployed in challenging environments for performing complex tasks such as multi-target tracking, search-and-rescue, and intrusion detection. Notwithstanding the computational limitations of individual robots, such systems rely on collaboration to sense and react to the environment. This paper formulates the generic target tracking problem as a time-varying optimization problem and puts forth an inexact online gradient descent method for solving it sequentially. The performance of the proposed algorithm is studied by characterizing its dynamic regret, a notion common to the online learning literature. Building upon the existing results, we provide improved regret rates that not only allow non-strongly convex costs but also explicating the role of the cumulative gradient error. Two distinct classes of problems are considered: one in which the objective function adheres to a quadratic growth condition, and another where the objective function is convex but the variable belongs to a compact domain. For both cases, results are developed while allowing the error to be either adversarial or arising from a white noise process. Further, the generality of the proposed framework is demonstrated by developing online variants of existing stochastic gradient algorithms and interpreting them as special cases of the proposed inexact gradient method. The efficacy of the proposed inexact gradient framework is established on a multi-agent multi-target tracking problem, while its flexibility is exemplified by generating online movie recommendations {for Movielens $10$M dataset}.  

\end{abstract}
	\begin{IEEEkeywords}
	Time varying optimization, stochastic optimization, target tracking, gradient descent methods. 
\end{IEEEkeywords}
\vspace{0mm}
	\section{ Introduction}\label{intro}
Multi-agent systems involve teams of robots capable of accomplishing complex tasks in a coordinated manner \cite{dubowsky2005concept,stone2000multiagent,ota2006multi}. Fueled by advances in sensing and communications, multi-agent systems are increasingly being used for challenging tasks such as multi-target tracking \cite{derenick2009convex,tu2013influence}, planetary exploration and mapping \cite{Trafton}, search-and-rescue, and intrusion detection \cite{abbeel2007application}. Achieving such team-based goals requires the robotic platforms to not only sense and understand the environment, but also cooperate among themselves through judicious data exchange and fusion. Consequently, resource allocation and optimization becomes an important aspect  of the overall motion planning problem. Indeed, the recent trend is to formulate target tracking as a constrained convex optimization problem that must be solved at every time step \cite{zinkevich2003online,hall2015online,simonetto2017decentralized,mokhtari2016online,zhang2016improved} . 

Such time-varying optimization problems have their origins in the control theory literature, where they have been applied to path planning \cite{uav} and dynamic parameter tracking problems \cite{TV_LMS,zhou2011multirobot}. Given the limited computational and communications capabilities of the mobile robots, solving the full per-time instant optimization problem before taking the action may not necessarily be viable. Instead, recent works have advocated the use of simpler one-iteration algorithms such as the online interior point, prediction-correction, online ADMM, and gradient descent methods, that have been shown to approach the optimal asymptotically. Leveraging the tools from classical optimization theory, these dynamic optimization algorithms not only admit low-complexity distributed implementations, but are also amenable to analytical performance guarantees. 

Online machine learning represents a parallel but closely related development that has been widely applied to solve problems in Big Data\cite{perera2009clustering}. First proposed by \cite{zinkevich2003online}, the online convex optimization framework models the agents as learners and targets as adversaries. Within this sequential learning paradigm, the learner performs an action and the adversary reveals a corresponding loss function at each time step. The eventual goal of the learner is to minimize the cumulative loss. Recent year have witnessed the development of theoretical guarantees in form of \emph{dynamic regret}, where the performance of the learner is measured against that of an adaptive and time-varying adversary \cite{besbes2015non,mokhtari2016online,shahrampour2016distributed}. 

This paper studies the multi-agent multi-target tracking problem from the lens of online convex optimization. Prompted by the noisy and possibly incorrect target position information available to the agents, we put forth the inexact online gradient descent (IOGD) algorithm. The key theoretical contribution is the development of improved dynamic regret bounds for non-strongly convex objective functions. Improving the existing results for general convex problems, it is shown that the dynamic regret is bounded by the \emph{path length} of the target. Different from the existing literature, the gradient errors are not required to be unbiased, independent, or identically distributed, and can even be adversarial. The cumulative impact of such errors on the dynamic regret is explicitly studied and characterized. Regret bounds are developed for two classes of problems, namely the ones where the objective function follows the quadratic growth property ($\pa$), and the ones where the function is convex but with a compact domain ($\pb$). Interestingly, it is shown that for $\pa$ with linearly growing path lengths and cumulative errors, target tracking is still possible, albeit with non-vanishing steady state error. Additionally, we provide further motivation for the IOGD algorithm, by developing online variants of the incremental algorithm from \cite{IGM_hy} and the proximal gradient method \cite{schmidt2011convergence}, that may still be analyzed within the current framework. Finally, the flexibility of the IOGD algorithm is demonstrated by applying it to the multi-agent multi-target tracking problem from \cite{derenick2009optimal} and the online matrix completion problem \cite{dhanjal2014online}. 

The rest of the paper is organized as follows. We begin with a brief review of literature on the areas of target tracking and online optimization. Sec.\ref{Prob_for} formulates the problem at hand and provides a succinct comparison of the existing results and those developed in this paper. The required regularity conditions and assumptions, as well as the proposed IOGD algorithm are stated in Sec.\ref{reg_1} and Sec.\ref{seciii}. Regret bounds for the two classes of problems are developed in Sec.\ref{rbp1} and Sec.\ref{rbp2}. In the general case when the path lengths or the cumulative errors are not sublinear, bounds on the steady state tracking error are developed in Sec.\ref{track}. The IOGD algorithm is also shown to subsume the online variants of two of the existing gradient algorithms in Sec.\ref{special}. Finally, numerical tests on real and synthetic data are provided in Sec.\ref{num}. 
\vspace{0mm}\subsection{Related work}
Inexact gradient methods have been widely used to solve a variety of optimization problems, especially in the context of machine learning \cite{IGM_p,IGM_rs, IGM_hy,lan2009convex}. Since calculating an approximate gradient is often cheap, recent years have witnessed the development of several variants, such as the incremental aggregated gradient method \cite{blatt2007convergent}, stochastic average gradient method \cite{roux2012stochastic,Nedic2009ApproximatePS}, stochastic variance reduced gradient method \cite{johnson2013accelerating}, SAGA \cite{defazio2014saga}. For static optimization problems, the inexact gradient methods are known to converge at a linear rate even for non-strongly convex objectives \cite{IGM_ppr}, \cite{luo1992linear}. \colb{It is remarked that the present work considers an online setting characterized by streaming or sequentially arriving data. Therefore, the performance of the proposed algorithm is measured against that of a slowly moving clairvoyant. In contrast, the inexact gradient methods (IGM) developed in \cite{IGM_ppr,IGM_p,IGM_rs,IGM_hy}are essentially offline algorithms that are meant to be applied when full data is available. The claims in \cite{IGM_ppr,IGM_p,IGM_rs,IGM_hy} therefore concern the asymptotic behavior of the iterates $\{\x_k\}$ and are quite different from those in the present paper. For instance, the dynamic regret bounds obtained here do not directly follow from the convergence results in these works.}

Time-varying optimization problems have classically been studied in the context of optimal control \cite{hours2016parametric}, target tracking \cite{dontchev2013euler,zavala2010real}, non-stationary optimization \cite{popkov2005gradient}, and parametric programming \cite{dontchev2009implicit}. First order gradient based methods have been advocated as efficient solvers for such problems \cite{zinkevich2003online, dyn_stoch, hall2015online, besbes2015non, jadbabaie2015online, mokhtari2016online,simonetto2017time,simonetto2017decentralized,yang2016tracking,zhang2016improved}, and will be discussed in detail in Sec. \ref{relreg}. When the dynamics of the target are partially known and the cost function is sufficiently smooth, it may be possible to use second order methods, such as those proposed in \cite{PC_y, PC_ppr, PC_ppr2, PC_ppr3,fazlyab2016prediction}. On the other hand, tracking is possible even for non-differentiable but strongly convex function using the subgradient and alternating directions method of multiplier methods \cite{simonetto2015non}. It is remarked that an underlying assumption in all these papers is that of strong convexity of the cost function. Such an assumption is quite strong and is not satisfied by several problems, such as least squares and logistic regression \cite{IGM_ppr}.

Dynamic regret for analyzing tracking problems was first introduced in \cite{zinkevich2003online}. As compared to the weaker notion of static regret, the idea here is to compare the performance of the tracker against that of an adaptive and time-varying adversary \cite{besbes2015non,mokhtari2016online,zhang2016improved}. Dynamic regret bounds for the gradient descent and related first order methods have been reported in \cite{zinkevich2003online,hall2015online,besbes2015non,jadbabaie2015online,mokhtari2016online,yang2016tracking,zhang2016improved}. An even stronger notion of offline regret has recently been introduced in \cite{chen2017online}. The present work builds upon the dynamic regret and steady state tracking results reported in \cite{mokhtari2016online,simonetto2015non,simonetto2017decentralized}. As compared to existing results, the bounds provided here are not only stronger but also require minimal assumptions on the cost function. 

\begin{table*}
	\centering
	\captionof{table}{Summary of related works  on time-varying optimization (cf. Sec.\ref{seciii})}	\label{table}
	\begin{tabular}{cccccc}
		\hline
		References &  Loss function & Inexact & Function & Regret rate\\
		\hline\\ 
		\cite{zinkevich2003online}  & Convex 			 	& No	& Deterministic  &	$\mathcal{O}{(\sqrt{K}(1+W_K))}$\\
		\cite{hall2015online} 		& Convex				& No	& Deterministic &	$\mathcal{O}{(\sqrt{K}(1+W_K))}$\\ 
		\cite{besbes2015non} 		& Convex 				& No	& Stochastic &	$\mathcal{O}{({K^{2/3}}(1+V_K)^{1/3})}$\\
		\cite{besbes2015non} 		& Strongly convex	 	& No	& Stochastic &	$\mathcal{O}{(\sqrt{K(1+V_K)})}$\\
		\cite{jadbabaie2015online} 	& Convex  			& No	& Deterministic &	$\mathcal{O}{(\sqrt{D_K+1}+\min\{\sqrt{(D_K+1)W_K},[(D_K+1)V_K K]^{1/3}\})}$\\
		\cite{mokhtari2016online} 	& Strongly convex 		& No	& Deterministic  &	$\mathcal{O}{(1+W_K)}$\\
	\cite{yang2016tracking} 	& Convex		& No	& Stochastic  &	${\mathcal{O}{(\sqrt{K W'_K})}}$\\
		\cite{zhang2016improved} 	& Convex+QG		& No	& Deterministic  &	${\mathcal{O}{(1+W''_K)}}$\\
		This work 					& Convex + QG ($\mathsf{P_1}$) 			 	& Yes	&  Deterministic/Stochastic (\textbf{A1})&	${\mathcal{O}{(1+E_K+W_K)}}$ \\
		This work 					& Convex ($\mathsf{P_2}$)  			 	& Yes	&  Deterministic/Stochastic (\textbf{A1}) &	${\mathcal{O}{(1+\sqrt{KE_K}+W_K)}}$ \\
		This work 					& Convex ($\mathsf{P_2}$)			 	& Yes	& Stochastic (\textbf{A2}) &	 ${\mathcal{O}{(1+E_K+W_K)}}$\\
		\hline
	\end{tabular}
	\vspace{0mm}
\end{table*}	
\textbf{Notations:} Scalars are denoted by letters in regular font, while vectors (matrices) are denoted by bold face (capital) letters. The $(i,j)$-th element of a matrix $\E$ is denoted by $[\E]_{ij}$. The all-one vector of size $n\times 1$ is represented by $\mathbf{1}_n$, while $\mathbf{I}_n$ denotes identity matrix of size $n \times n$. The notation $\norm{\cdot}$  represents the Euclidean norm. The Kronecker product operator is denoted by $\otimes$. The expectation operator is symbolized by $\mathbb{E}$. Finally, $\sigma_{\max}(\E)$ and $\sigma_{\min}(\E)$ represent the maximum and minimum non-zero singular values of matrix $\E$, respectively. 
\vspace{0mm}
\section{Problem Formulation}\label{Prob_for}
This section details the target tracking problem and formulates it within the online learning framework. An example is also provided to motivate the setting. 

\subsection{System Model}
We consider the general problem of tracking a time-varying parameter that evolves according to an unknown dynamic model. The general setting considered here subsumes the target tracking application of interest, where the parameter may represent the location(s) of the target(s) being pursued by one or more agents. As motivated in \cite{mokhtari2016online,simonetto2015non, shahrampour2016distributed}, the parameter at time $k$ can be written as the solution of the following (discrete) time-varying convex optimization problem 
\begin{align}
\x_k^\star &\in \arg \min_{\x \in \cX} f_k(\x) & k = 1, 2, \ldots \label{xkstar}
\end{align}
where $f_k$ is a smooth convex function and $\cX \subset \Rn^n$ is a convex set. The set notation in \eqref{xkstar} emphasizes the fact that in general, the minimizer of $f_k$ may not necessarily be unique. The parameter estimate at time $k$ is denoted by $\x_k$. The agents do not know the full functional form of $f_k$ but are only revealed an inexact version of the gradient $\tilde{\nabla} f(\x_k) := \nabla f_k(\x_k) + \e_k$ for some $\e_k \in \Rn^n$. The agents make use of these inexact gradients to improve their estimates of $\x_k^\star$ in an online manner.

This paper considers the problem from an online convex optimization perspective, viewing the agents as learners and targets as adversaries. Specifically, at time $k$, the online learner selects an action $\x_k \in \cX$ and incurs a cost $f_k(\x_k)$, where $f_k:\Rn^n \rightarrow \Rn$ are smooth convex functions selected by the adversary. In response to the agent's action, the adversary also reveals an inexact gradient $\tilde{\nabla} f_k(\x_k)$ to the learner. The goal of the learner is to minimize its cumulative loss $\sum_{k=1}^K f_k(\x_k)$ over $K$ time slots. Of particular interest is to characterize the so-called \emph{dynamic regret} of the learner, that measures the cumulative mismatch between the learner's action and the optimal action \cite{mokhtari2016online,besbes2015non,hall2015online}:
\begin{align}\label{reg}
\reg_K:=\sum\limits_{k=1}^{K}(f_k(\x_k)-f_k(\x_k^\star))
\end{align}
where $\x_k^\star$ is as defined in \eqref{xkstar}. In order for the tracking to be successful, the dynamic regret must be upper bounded by a sublinear function of $K$. 
\vspace{0mm}
\subsection{Parameter variations and error bounds}\label{reg_1}
It is well known that a sublinear dynamic regret may not always be achievable, e.g., if the parameter variations or the gradient errors are too large \cite{besbes2015non,zhang2016improved}. The goal of the current paper will therefore be to bound the dynamic regret using functions of the cumulative parameter variations and errors. For the target tracking setting at hand, it makes sense to characterize the parameter variations using the path length, defined as
\begin{align}\label{pl}
W_K:=\sum\limits_{k=2}^{K}\norm{\x_{k}^\star-\x_{k-1}^\star}
\end{align} 
for some sequence of parameter values $\{\x_k^\star\}_{k\geq 1}$. More generally there exist a class of related complexity measures that can be used to characterize the parameter variations \cite{mokhtari2016online}. Examples include the functional variation $V_K:=\sum\limits_{k=2}^{K}\sup_{\x\in\cX}\abs{f_k(\x)-f_{k-1}(\x)}$ and the gradient variation $D_K:=\sum\limits_{k=2}^{K}\sup_{\x\in\cX}\norm{\nabla f_k(\x)-\nabla f_{k-1}(\x)}^2$. 

The gradient errors $\e_k$ can be modeled either as being deterministic with bounded norms or stochastic with bounded variance. Deterministic errors may be of interest in adversarial settings while stochastic errors are useful for modeling communication and computational noise. In order to unify the subsequent development, a generic stochastic error model is considered that subsumes that deterministic case. Let $\Fk$ denote the sigma field generated by the random sequence $\{\e_\tau\}_{\tau = 1}^{k-1}$. The following assumption bounds the second moment of the error. 

\noindent\emph{\noindent\textbf{A1. Error bound}: The stochastic sequence $\{\e_k\}_{k=1}^K$ adheres to the following bound on the second moment:
\begin{align}
\Ex{\norm{\e_k}^2~|~\Fk}&\leq \varepsilon_k^2+\nu^2\norm{\nabla f_k(\x_k)}^2 \label{a1}
 \end{align}
with $\varepsilon_k\leq \varepsilon < \infty$ for all $k \geq 1$, and $\nu \geq 0$ and $\varepsilon \geq 0$ are constants}.

When the errors are deterministic, \eqref{a1} is equivalent to a worst-case bound on the error norm. In the general case, the goal will be to establish the sublinearity of the expected dynamic regret $\Ex{\reg_K}$. The specific form of the bounds in \textbf{(A1)} is inspired from \cite{bertsekas2000gradient, rosasco2014convergence} and allows for errors that are proportional to the gradient norm in addition to an additive term. It is remarked that from Jensen's inequality, \eqref{a1} implies that $\Ex{\norm{\e_k}\mid \Fk}\leq \varepsilon_k+\nu\norm{\nabla f_k(\x_k)}$. The required dynamic regret bounds will be developed in terms of the path length $W_K$ and the cumulative error bound $E_K := \sum_{k=1}^K \varepsilon_k$.	

A particular case of interest is when the gradient errors constitute a white noise process as specified in the following assumption.

\noindent\emph{\noindent\textbf{A2. White noise}: The zero-mean error sequence $\{\e_k\}_{k=1}^K$ is independent identically distributed, i.e.,
\begin{align}\label{a2}
\Ex{\e_k \mid \Fk} = \mathbf{0}.
\end{align}}
Assumption (\textbf{A2}) may be applicable, for instance, when the gradient errors arise from communication errors. The requirement in \eqref{a2} is more restrictive than that in \eqref{a1}, but also results in improved regret bounds.

The dynamic regret bounds will only be meaningful when the quantities $W_K$ and $E_K$ are sublinear in $K$. In the context of target tracking, a sublinear $W_K$ is obtained, for instance, if the speed of the target goes to zero as $K$ increases. Such situations are commonly encountered for targets with finite total operational energy. Likewise, $E_K$ may be sublinear in scenarios where the gradient error goes to zero with $K$. Again, such situations may arise if the noise can be filtered out with increasingly higher accuracy. Finally, for the case when the $W_K$ and $E_K$ are linear functions of $K$, it makes sense to instead characterize the tracking performance $\norm{\x_k-\x_k^\star}$ of the proposed algorithm. Different from regret, the goal here is to simply ensure that the asymptotic tracking performance is not too large; see Sec. \ref{track}. 
\vspace{0mm}
\subsection{Related work and problem statement} \label{relreg}
Before detailing the problem statement, the existing literature on dynamic regret bounds is briefly reviewed. In \cite{zinkevich2003online}, an online gradient descent algorithm is considered and shown to achieve $\mathcal{O}(\sqrt{K}(1+W_K))$ regret bound under diminishing step size. Similar regret bounds are shown to exist for the online mirror descent algorithm in \cite{hall2015online}. The online gradient descent with stochastic error is considered in \cite{besbes2015non} and  regret bounds for both convex and strongly convex case are provided as $\mathcal{O}{({K^{2/3}}(1+V_K)^{1/3})}$ and $\mathcal{O}{(\sqrt{K(1+V_K)})}$, respectively. It is remarked that these results cannot be practically applied to tracking problems since they assume a priori knowledge of the functional variation $V_K$. A regret bound in terms of all the three measures $W_K$, $V_K$ and $D_K$ is obtained in \cite{jadbabaie2015online} for optimistic mirror descent algorithm. In succession to these works, authors in \cite{mokhtari2016online,zhang2016improved} achieved the $\mathcal{O}(1+W_K)$ dynamic regret bound for online gradient descent algorithm under constant step size and strongly convex functions. \colb{The case of non-strongly convex functions has been recently considered in \cite{zhang2016improved} and a regret bound of $\mathcal{O}(1+W_K'')$ is obtained for deterministic time-varying problems. The present results generalize those in \cite{zhang2016improved} to stochastic problems with adversarial gradient errors, unbounded domain $\mathcal{X}$, and a stronger path length definition \eqref{pl}. Finally, the case of noisy gradient is also considered in \cite{yang2016tracking} for general convex functions and dynamic regret bounds are obtained while assuming the knowledge of $K$. The results in the present work improve upon those in \cite{yang2016tracking}, in addition to allowing gradient errors that are not necessarily independent or identically distributed (cf. (\textbf{A1}).  Moreover, the impact of gradient errors on the regret is explicitly specified, allowing us to develop IOGD variants that cannot be analyzed as special cases of the noisy OGD proposed in \cite{yang2016tracking}.} Table \ref{table} summarizes the regret bounds obtained in the existing literature. Next, we briefly remark on the path length definitions used in \cite{yang2016tracking,zhang2016improved}. 

\begin{rem}\label{remark1}
The path length definition used in \eqref{pl} applies to an arbitrary sequence of true parameter values $\{\x_k^\star\}$ and does not depend on $\cX_k^\star$. Consequently, the path length definition in \eqref{pl} is stronger than those used in \cite{yang2016tracking, zhang2016improved}. In particular, the definitions in \cite{yang2016tracking,zhang2016improved} take the following form
\begin{align}
W'_K :=& \max_{\col{\{u_k \in \Xk\}_{k=2}^{K}}}\sum\limits_{k=2}^{K}\norm{\u_{k}-\u_{k-1}}\label{rem_first}\\
W''_K :=& \sum\limits_{k=2}^{K}\max_{{\u\in\mathcal{X}}}\norm{P_{k}(\u)-P_{k-1}(\u)},\label{rem_second}
\end{align}
respectively, where $P_k(\u):=\arg\min_{\y\in\mathcal{X}_k^\star}\norm{\y-\u}$. Observe that in both cases, the path length could become unbounded if the set $\cX_k^\star$ is not compact, as is the case in $\pa$. \colb{For example, the path lengths $W'_K$ and $W''_K$ turn out to be  infinite for the least squares example discussed in Sec. \ref{time_varying_estimation}}. The use of an arbitrary trajectory in \eqref{pl} allows us handle such unbounded sets while also ensuring that $W_K \leq W'_K$ and $W_K \leq W''_K$ for problems where $\cX_k^\star$ is compact. 
\end{rem}

The present work considers the online inexact gradient descent method and develops various regret and tracking bounds for the case when the cost function $f_k$ is not strongly convex. The bounds developed here are efficient, and match the existing \colb{bounds in the cases} when the gradient is exact. Keeping the tracking application in mind, the algorithms developed here will not require any prior knowledge of $K$. To this end, we consider two broad classes of non-strongly convex functions. 

\noindent \textbf{$\pa$. Special structure:} \emph{The function $f_k$ takes the form $f_k(\x) = g_k(\E\x)$ where $\E$ is an arbitrary $m \times n$ matrix and the function $g_k$ is strongly convex with parameter $\mu >0$. }

The functional form required in {$\pa$} is significantly more general than the strong convexity condition for $f_k$. Since $g_k$ is strongly convex, there exists a unique $\u^\star  = \arg\min_{\y \in \cX} g_k(\y)$. Consequently, the set of minimizers of $f_k$ is given by $\Xk:=\{\x_k^\star \mid \E\x_k^\star = \u^\star\}$. The key property of $f_k$ that will subsequently be utilized is the so-called quadratic growth condition\vspace{-0mm}
\begin{align}\label{eq:QG}
f_k(\x)-f^\star_k &\geq \frac{\mu \sigma_{\min}(\E)}{2} \norm{\x-P_k(\x)}^2
\end{align}
where \colb{$f^\star_k:=f_k(\x_k^\star)$ for $\x_k^\star\in\mathcal{X}_k^\star$}, $ P_k(\x):=\arg\min_{\y\in\mathcal{X}_k^\star}\norm{\y-\x}$ and $\sigma_{\min}(\E)$ denotes the smallest non-zero singular value of $\E$. The quadratic growth condition for convex problems is weaker than the standard assumption of strong convexity \cite{karimi2016linear}. For instance, $\pa$ is applicable to  common machine learning and signal processing problems such as least squares, support vector machines, and logistic regression, none of which are strongly convex \cite{IGM_ppr}. 

\noindent \textbf{$\pb$. Bounded domain: } \emph{For each $k\geq 1$, $f_k$ is convex with the set of minimizers denoted by $\Xk$ and the set $\cX$ is compact with diameter $R < \infty$. }

 It is remarked that the bounds developed for $\pb$ are directly proportional to the diameter $R$, and may therefore be quite loose. Before concluding, a remark on the possible generalization of $\pa$ is due. 
\begin{rem}
The results developed here can be extended to the generalized version of $\pa$ where $f_k(\x) = g_k(\h(\x))$ satisfies the QG condition, $\h$ is a smooth function, and $g_k$ is strictly convex.  Alternatively, it is also possible to consider convex QG functions $f_k$ that have a unique minimum. The extensions can be incorporated via minor modifications in the proofs, but are not pursued here in order to keep the exposition succinct.
\end{rem}\vspace{-0mm}
\section{Proposed algorithm and assumptions}\label{seciii}
The online gradient descent algorithm has been widely used to solve online learning problems owing to its flexibility and simplicity \cite{zinkevich2003online,hall2015online,besbes2015non,jadbabaie2015online,mokhtari2016online}. This work considers the inexact online gradient descent (IOGD) method that takes the form:
\begin{equation}
\label{IGM}
\x_{k+1} = \col{\mathcal{P}_{\mathcal{X}}}[\x_{k} - \alpha (\nabla f_k(\x_k) + \e_k)]
\end{equation}
where $\col{\mathcal{P}_{\mathcal{X}}}(\cdot)$ denotes the projection onto the set $\cX$ \colb{and is defined as $\mathcal{P}_\mathcal{X}(\u):=\arg\min_{\y\in\mathcal{X}}\norm{\y-\u}$}. The IOGD method has also been applied to static and online problems \cite{IGM_ppr}. The IOGD method is also closely related to the incremental and variance reduced variants of the gradient descent algorithm. The full algorithm is summarized in Algorithm \ref{algo_1}. 

	  \begin{algorithm}
	  	\caption{IOGD: Inexact Online Gradient Descent}\label{algo_1}
	  	\begin{algorithmic}[1]
	  		\STATE {\textbf Initialize} $\x_{1}$	  		
	  		\STATE {\textbf{for} $k=1$, $2$, $\ldots$} \textbf{do} 	  		
	  		\STATE \ \  \colb{\textbf{Perform} action $\x_k$}
	  		\STATE \ \ \textbf{Observe} inexact gradient $\nabla f_k(\x_k) + \e_k$ at $\x_k$  		
	  		\STATE \ \ \colb{\textbf{Compute} next action }\colb{$\x_{k+1} =\mathcal{P}_{\mathcal{X}}[\x_{k} - \alpha (\nabla f_k(\x_k) + \e_k)]$ }\vspace{-5mm}		
	  		\STATE \textbf{end for }
	  	\end{algorithmic}
	  		  	\label{algo1}
	  \end{algorithm}

In addition to Assumptions (\textbf{A1}) and (\textbf{A2}) stated in Sec. \ref{Prob_for}, the subsequent analysis will also require the following regularity conditions. 

\noindent\emph{\noindent\textbf{A3. Lipschitz continuity}: The function $\nabla f_k$ is Lipschitz continuous with parameter $L$:
\begin{align}
\norm{\nabla f_k(\u) - \nabla f_k(\v)}&\leq L\norm{\u-\v}  \label{Lip_g12}
\end{align}
for all $k \geq 1$ and $\u$, $\v \in \cX$.}

\noindent\emph{\noindent\textbf{A4. Vanishing gradient}: The optimum $\x_k^\star$ lies in the relative interior of the set $\cX$, that is, $\nabla f_k(\x_k^\star) = 0$ for all $\x_k^\star \in \Xk$.}

\noindent\emph{\textbf{A5. Bounded Variation}}: For a given $\x_k^\star$, there exists some $\sigma>0$ such that $\norm{\x_{k+1}^\star-\x_{k}^\star}\leq\sigma$ for all $k\in\mathbb{N}$.

Of these, both (\textbf{A3}) and (\textbf{A4}) are standard and apply to large class of online learning problems. Likewise, the requirement in (\textbf{A5}) imposes a limit on the maximum velocity of the target and is therefore applicable to most target tracking problems. The bounded variation condition is also satisfied, for instance, if $W_K$ is sublinear or linear and the target motion is not too `jumpy'. 

The subsequent two sections develop the dynamic regret bounds stated in Table \ref{table}. As stated earlier, these bounds are meaningful when the variations $W_K$ and error bounds $E_K$ are sublinear in $K$. Development of the tracking error bounds for linear $W_K$ and $E_K$ is deferred to Sec. \ref{track}.

\section{Regret bounds for $\pa$} \label{rbp1}
This section develops the dynamic regret bounds for $\pa$ for the case when $W_K$ and $E_K$ are sublinear in $K$. As discussed earlier, a simple example is that of a target that eventually stops, either upon exhausting its energy or upon reaching its destination. For such targets, there exists $K_0 < \infty$ such that $\x_{k+1}^\star=\x_k^\star$ for all $k\geq K_0$, making $W_K$ constant for $K \geq K_0$. Likewise, the  term $\varepsilon_k$ in the gradient error bound may decrease over time resulting in a sublinear $E_K$.  Interestingly however, the gradient error is allowed to be proportional to the gradient norm, and the factor $\nu$ need not be diminishing with $k$.  

The section proceeds by developing regret bounds for the general case in (\textbf{A1}), while the special case when the gradient errors constitute a white noise process as in (\textbf{A2}) is treated towards the end. \colb{The mathematical analysis for these two cases is unified using an indicator variable $1_d$, that takes the value $1$ when (\textbf{A1}) is in effect and zero when both (\textbf{A1}) and (\textbf{A2}) are in effect (cf. Appendix \ref{proof_lem1}). } The overall proof involves four key steps, namely, (a) bounding the quantity $\Ex{\norm{\x_{k+1} - P_k(\x_{k})}}$ in terms of the gradient error and $\dist{\x_k,\Xk}:=\norm{\x_k-P_k(\x_k)}$; (b) bounding the cumulative sum of $\dist{\x_k,\Xk}$; (c) establishing a bound on the average gradient norm $\Ex{\norm{\nabla f_k(\x_k)}}$; and finally, (d) using the gradient norm bound to obtain the required bound on $\Ex{\reg_K}$. {} The following intermediate \colb{lemma} develop the bounds required in the first three steps. The proofs of these results are deferred to Appendix \ref{proof_lem1}. 

\begin{lem}\label{lema1}
Under \textbf{(A1)}, and \textbf{(A3)}-\textbf{(A5)}, the sequence  $\{\x_k\}$ for all $k\in\mathbb{N}$ generated by IOGD algorithm satisfies the following bounds for any sequence $\{\x_k^\star\}$ \colb{with sublinear path length $W_K$}:
	\begin{align}
	\Ex{\norm{\x_{k+1} - P_k(\x_k)}} &\leq {\ell}\Ex{\dist{\x_k,\cX_k^\star}}+{\zeta} \varepsilon_k\label{contraction} \\
	\sum_{k=1}^K\Ex{\dist{\x_k,\cX_k^\star}} &\leq \frac{\norm{\x_1-\x_1^\star} + W_K + {\zeta} E_K}{\chi-{\ell}} 	\label{sumdist} 
	\\
	\Ex{\norm{\nabla f_k(\x_k)}} &\leq G\!:=\!L\frac{{\zeta} \varepsilon \!+\! \sigma}{\chi} +L\dist{\x_1,\cX_1^\star}\label{gradient_bound}
	\end{align}
	where, \colb{${\ell^2}$ is obtained from \eqref{ell_definition} using $1_d=1$ as 
	\begin{align}
	\ell^2&:= 1-\mu\alpha(1-\alpha L(1+\nu)^2)+2\nu \alpha L  \nonumber
	\end{align}
	{and $\zeta$ is obtained from \eqref{zeta-updaet} using $1_d=1$ as}
	\begin{align}
	{\zeta}&:= \frac{\alpha (1+\alpha L)}{{\ell}}.\vspace{-0mm}
	\end{align}} 
	\colb{The step size $\alpha$ is chosen such that 
	\begin{align}\label{alpha_range}
	\frac{\mu-2\nu L-\varpi}{2\mu L(1+\nu)^2} < \alpha < \frac{\min\{\mu-2\nu L+\varpi,2\mu\}}{2\mu L(1+\nu)^2}
	\end{align}
	where $	\varpi^2 = (\mu-2\nu L)^2-4\mu L(1+\nu)^2(1-\chi^2)$. } 
\end{lem}

It is remarked that an appropriate value of $\alpha$ satisfying the conditions required in Lemma \ref{lema1} always exists if $\mu - 2\nu L > 2(1+\nu)\sqrt{\mu L (1-\chi^2)}$. For example, in the case when $\E$ is an identity matrix, the condition becomes $\mu > 2\nu L$, while in the case when $\nu = 0$, the condition translates to $\mu > 4L(1-\chi^2)$. \colb{Interestingly, the selection of $\alpha$ does not require the prior information about $K$. As also detailed in \cite{mokhtari2016online,simonetto2015class}, target tracking algorithms are generally required to run indefinitely or until the some stopping criteria such as the `distance to target is sufficiently small' is met. In either case, the number of iterations $K$ is generally not known or fixed in advance.} 

Among the bounds developed in Lemma \ref{lema1}, the gradient bound in \eqref{gradient_bound} is particularly interesting. Different from the vast majority of literature on online convex optimization, the gradient bound allows us to forgo the assumption on the boundedness of the set $\cX$ \cite{zinkevich2003online,hall2015online,besbes2015non}. Further, such a result is useful even when the set $\cX$ is compact but has a large diameter $R$, for the dynamic regret results presented here will not depend on $R$. Having stated the intermediate lemma, the main result of this section can finally be stated as the following theorem.

\begin{thm}\label{lem:deter_err}
Under (\textbf{A1}), (\textbf{A3})-(\textbf{A5}), \colb{and $\alpha$ satisfying $\eqref{alpha_range}$}, the IOGD iterates result in the following dynamic regret bound
\begin{align}\label{main_theorem}
\Ex{\reg_K} \leq \mathcal{O}(1+E_K+W_K).
\end{align}
\end{thm}

\begin{IEEEproof}
The proof follows from the first order convexity condition for $f_k$ and the use of the bounds in \eqref{sumdist} and \eqref{gradient_bound}. Given a sequence $\{\x_k^\star\}$ \colb{with sublinear path length $W_K$}, it holds that
\begin{align}
\sum_{k=1}^K (f_k(\x_k) - f_k(\x_k^\star)) &= \sum_{k=1}^K (f_k(\x_k) - f_k(P_k(\x_k))) \nonumber\\
&\leq \sum_{k=1}^K \ip{\nabla f_k(\x_k), \x_k-P_k(\x_k)} \\
&\leq \sum_{k=1}^K \norm{\nabla f_k(\x_k)}\norm{\x_k\!-\!P_k(\x_k)}
\end{align}
where we have used the first order convexity condition and the Cauchy-Schwarz inequality. Taking expectation and using the bounds in \eqref{sumdist} and \eqref{gradient_bound}, it follows that
\begin{align}
\Ex{\reg_K} &\leq  \left(L\frac{{\zeta} \varepsilon + \sigma}{\chi-{\ell}} +L\dist{\x_1,\cX_1^\star}\right)\times\nonumber\\
&\hspace{1cm}\left(\frac{\norm{\x_1-\x_1^\star} + W_K + {\zeta} E_K}{\chi-\ell}\right) \\
&\leq \mathcal{O}(1+W_K+E_K)
\end{align}
which is the required result. 
\end{IEEEproof}	  		

The behavior of regret bound for the IOGD algorithm is governed by the target trajectory $\{\x_k^\star\}$ and the error sequence $\e_k$. The results provided here subsume the OGD results in \cite{mokhtari2016online} where $f_k$ is strongly convex and the exact gradient is available, i.e., $E_K = 0$. Interestingly, this \colb{theorem} establishes that an $\mathcal{O}(1+W_K)$ is still obtainable for $\pa$ with inexact gradients as long as $E_K \leq \mathcal{O}(W_K)$.

Finally, we provide the dynamic regret bounds for the case when the gradient errors follow a white noise process as in (\textbf{A2}). 				

\begin{cor} \label{p1iidcor}
Under (\textbf{A1})-(\textbf{A5}) \colb{and $\alpha$ satisfying \eqref{alpha_range}}, the dynamic regret for the IOGD iterates is bounded as follows
\begin{align}
\Ex{\reg_K} &\leq \left(L\frac{\alpha \varepsilon + \sigma}{\chi-{\ell}} \right)\left(\frac{\norm{\x_1-\x_1^\star} + W_K + \alpha E_K}{\chi-{\ell}}\right) \nonumber\\
&=\mathcal{O}(1+W_K+E_K)
\end{align}
where \colb{$\ell$ is obtained from \eqref{ell_definition} with $1_d=0$ given by $\ell:=\sqrt{1-\mu\alpha(1-\alpha L(1+\nu^2))}$ and $\ell<\chi$}. \colb{The value of $\alpha$ is chosen such that 
	\begin{align}
	\frac{\mu-\varpi}{2\mu L(1+\nu^2)} < \alpha < \frac{\min\{\mu+\varpi,2\mu\}}{2\mu L(1+\nu^2)}\nonumber
	\end{align}
	where $	\varpi^2 = (\mu)^2-4\mu L(1+\nu^2)(1-\chi^2)$. }
\end{cor}

\begin{IEEEproof}
As shown in Appendix \ref{proof_lem1}, the intermediate results for (\textbf{A1})-(\textbf{A2}) take the form
\begin{align}
	\Ex{\norm{\x_{k+1} - P_k(\x_k)}} &\leq {\ell}\Ex{\dist{\x_k,\cX_k^\star}}+\alpha \varepsilon_k\label{contraction-iid} \\
	\sum_{k=1}^K\Ex{\dist{\x_k,\cX_k^\star}} &\leq \frac{\norm{\x_1-\x_1^\star} + W_K + \alpha E_K}{\chi-{\ell}} 	\label{sumdist-iid} \\
	\Ex{\norm{\nabla f_k(\x_k)}} &\leq L\frac{\alpha \varepsilon + \sigma}{\chi-lr{\ell}} +L\dist{\x_1,\cX_1^\star}\label{gradient_bound-iid}
	\end{align}
	where ${\ell}$ is as defined in Corollary \ref{p1iidcor}. The required dynamic regret bounds can be obtained by proceeding along the lines of proof of Theorem \ref{lem:deter_err}. 	
\end{IEEEproof}
It is remarked that valid values of $\alpha$ exist if $\mu > 4L(1-\chi^2)(1+\nu^2)$. It can be observed that the additional assumption in (\textbf{A2}) does not provide any improvement in the asymptotic regret rate, but only in the associated constants.

\section{Regret Bound for $\mathsf{P_2}$} 	\label{rbp2}
Having established the regret rates for $\pa$ with a special QG structure, we now move to the more general problem in $\pb$. However, for this case, it would not be possible to develop the results for generic stochastic errors following (\textbf{A1}). Instead, the \colb{results will be} presented here for two specific scenarios, namely (a) gradient errors following (\textbf{A1}) with $\nu = 0$, and (b) gradient errors following (\textbf{A1})-(\textbf{A2}) but possibly non-zero value of $\nu$.  We begin with stating the following intermediate lemma whose proof is deferred to Appendix \ref{proof_deter_general}. 

\begin{lem}\label{lam_gen_deter}
Under (\textbf{A1}) with $\nu = 0$, (\textbf{A3})-(\textbf{A4}), and for a sequence $\{\x_k^\star\}$ \colb{with sublinear path length $W_K$},  it holds for $\pb$ that the IOGD iterates satisfy
\begin{align}
&\Ex{\norm{\x_{k+1}\!-\!\x_{k}^{\star}}} \leq \Ex{\norm{\x_k\!-\!\x_k^\star}} \!+\! \frac{{\xi}}{R}(\colb{\Ex{f_k(\x_k)\!-\!f_k(\x_k^\star)}}) + s_k\nonumber
\end{align}
where \colb{$\xi:=2\alpha(1-2\alpha L)$ and $s_k:=\sqrt{2\alpha^2\varepsilon_k^2 +2\alpha\varepsilon_kR}$ (cf. \eqref{lemm2_ini1} and \eqref{lemm2_ini2})  under $1_d=1$. The step size $\alpha$ is chosen as $0<\alpha<1/2L$}.
\end{lem}

Lemma \ref{lam_gen_deter} leads directly to the required dynamic regret bounds for $\pb$ under (\textbf{A1}) with $\nu = 0$. 

\begin{thm}\label{gen_det}
Under (\textbf{A1}) with $\nu = 0$, (\textbf{A3})-(\textbf{A4}), and for a sequence $\{\x_k^\star\}$ \colb{with sublinear path length $W_K$}, it holds for $\pb$ that the IOGD iterates adhere to the following dynamic regret rate
\begin{align}
\Ex{\reg_K} &\leq \mathcal{O}(1+\sqrt{KE_K}+W_K).
\end{align}
\colb{with step size selected as $0<\alpha<1/2L$.}
\end{thm}
\begin{IEEEproof}
Using triangle inequality and result of lemma \ref{lam_gen_deter} for $k \geq 1$, we have that 
\begin{align}
\Ex{\norm{\x_{k+1}-\x_{k+1}^\star}} &\leq \Ex{\norm{\x_{k+1} - \x_k^\star}}  + \norm{\x_{k+1}^\star - \x_k^\star}\nonumber\\
 	& \leq \Ex{\norm{\x_k - \x_k^\star}}-  \frac{{\xi}}{R}\Ex{f_k(\x_k)\!\!- \!\!f_k(\x_k^\star)}\nonumber\\ 
 	&\ \ \ \ \  + s_k + \norm{\x_{k+1}^\star - \x_k^\star}. \label{decrement}
 	\end{align}
It is now possible to apply the relationship 	in \eqref{decrement} recursively for $k = 1,\ldots, K$, so as to obtain
\begin{align}
\Ex{\norm{\x_{K+1}\!-\!\x_{K+1}^\star}} \!\leq &\norm{\x_1 \!-\! \x_1^\star}- \frac{{\xi}}{R}\!\sum_{k=1}^K\! \Ex{f_k(\x_k)\!-\! f_k(\x_k^\star)}  \nonumber\\
& + \sum_k s_k + W_K.
\end{align}
Since the left hand side is positive by definition, it follows that 
\begin{align}
\!\!\!\!\frac{{\xi}}{R}\Ex{\sum_{k=1}^K (f_k(\x_k)- f_k(\x_k^\star))} \leq \norm{\x_1-\x_1^\star} + S_k + W_K \label{regsk}
\end{align}
where we have that 
\begin{align}
S_K &= \sum_{k=1}^K \sqrt{2\alpha^2\varepsilon_k^2 + 2\alpha\varepsilon_kR} \\
&\leq  \sum_{k=1}^K \sqrt{2\alpha R\varepsilon_k} + \sqrt{2}\alpha\varepsilon_k \\
&\leq \sqrt{2\alpha R}\left(K\sum_{k=1}^K \varepsilon_k\right)^{1/2} + \sqrt{2}\alpha E_K \\
& = \sqrt{2\alpha R K E_K} + \sqrt{2}\alpha E_K.\label{skbound} 
\end{align}
Substituting the result in \eqref{skbound} into \eqref{regsk} yields the required result
\begin{align}
\Ex{\reg_K} &\leq \frac{R}{{\xi}}\left(\norm{\x_1\!-\!\x_1^\star} \!+\! \sqrt{2\alpha R KE_K} \!+\! \sqrt{2}\alpha E_K \!+\! W_K\!\!\right) \nonumber\\
&\leq \mathcal{O}(1+\sqrt{KE_K} + W_K)\label{here}
\end{align}
where the inequality in \eqref{here} follows since $\sqrt{KE_K} > E_K$ whenever $E_K$ is sublinear. 
\end{IEEEproof}
The results can be improved for the case when the gradient errors follow a white noise process. Specifically, the following corollary holds as shown in Appendix \ref{proof_deter_general}. 
\begin{cor}\label{cor_gen_stoc}
Under (\textbf{A1})-(\textbf{A4}), and for a sequence $\{\x_k^\star\}$ \colb{with sublinear path length $W_K$}, it holds for $\pb$ that the IOGD iterates satisfy
\begin{align}
&\Ex{\norm{\x_{k+1}\!-\!\x_{k}^{\star}}} \!\leq\! \Ex{\norm{\x_k\!-\!\x_k^\star}}\! +\! \frac{{\xi}}{R}(\colb{\Ex{f_k(\x_k)\!-\!f_k(\x_k^\star)}}) \!+\! {\alpha\varepsilon_k} \label{p2cor}
\end{align}
\colb{where ${\xi}:=2\alpha(1-\alpha L(1+\nu^2 ))$ since $1_d=0$ (cf. \eqref{lemm2_ini1}).} For this case, the dynamic regret is bounded as
\begin{align}
\Ex{\reg_K} \leq \mathcal{O}(1+E_K+W_K).
\end{align}
\colb{for $0 < \alpha<\frac{1}{L(1+\nu^2)}$}.
\end{cor}
\begin{IEEEproof}
The proof of the bound in \eqref{p2cor} is provided in Appendix \ref{proof_deter_general}. The required regret bound follows along the lines of the proof of Theorem \ref{gen_det}. \end{IEEEproof}
\vspace{0mm}
\section{ Tracking Performance}\label{track}
Departing from the regret analysis pursued in Sections \ref{rbp1} and \ref{rbp2}, we now consider the asymptotic tracking performance for the scenarios when $W_K$ or $E_K$ are not sublinear in $K$. As a practical example, consider a target that continues to move with a constant velocity in an adversarial manner. Likewise, we may consider the scenario when errors are unpredictably random with constant variance. More generally, this section will only assume that the variations and the gradient errors terms are bounded as follows

\begin{align}\label{boundedekwk}
\norm{\x_{k+1}^\star-\x_k^\star} &\leq \sigma \ \ \ \text{and} \ \ \ 
\varepsilon_k \leq \varepsilon
\end{align}
where realistically, the terms $\sigma$ and $\varepsilon$ should not be too large. In both cases, it may not be possible for the agents to reduce their distance from the target beyond a certain value \cite{simonetto2015class,simonetto2015non}.

\colb{The goal of this section will be to characterize the tracking performance $\Ex{\norm{\x_k-\x_k^\star}}$}. Towards this end, only the problem $\pa$ will be considered. As before, we begin with analyzing the general case when (\textbf{A1}) is in effect, and specialize the results later to the case when both (\textbf{A1})-(\textbf{A2}) are in effect. The main result is provided in the subsequent lemma, whose proof is provided in Appendix \ref{proof_lem1}. 

\begin{lem}\label{lem:tracking}
Under (\textbf{A1}) and (\textbf{A3})-(\textbf{A5}), and sequences $\{\x_k^\star, \e_k\}$ satisfying \eqref{boundedekwk}, the IOGD iterates for $\pa$ adhere to the following bound
\begin{align}\label{tracking_per0}
\Ex{\dist{\x_{k+1},\cX_{k+1}^\star}}\leq &\left({{\ell}}/{\chi}\right)^k\dist{\x_1,\cX_1^\star}  \nonumber
\\
&+ \left[{\frac{1-\left({{\ell}}/{\chi}\right)^k}{1-\left({\ell}/{\chi}\right)}}\right]\frac{\sigma +{\zeta} \varepsilon}{\chi}.
\end{align}	
for all $k \geq 1$. \colb{Here, the constants $\chi$, ${\ell}$, and ${\zeta}$ are as defined in Lemma \ref{lema1} {since $1_d=1$ for this case } {and $\alpha$ is chosen as in \eqref{alpha_range}}}.
\end{lem}

Lemma \ref{lem:tracking} establishes that the tracking error decreases exponentially with $k$, and is ultimately bounded by a steady state value of $\frac{\sigma +{\zeta} \varepsilon}{\chi-{\ell}}$. Alternatively, if it holds that $\Delta_1:=\dist{\x_1,\cX_1^\star} > \epsilon$, then the number of iterations required for the tracking error is to satisfy
\begin{align}\label{tracking_per1}
&\Ex{\dist{\x_{k+1},\cX_{k+1}^\star}} \leq \epsilon + \frac{\sigma +{\zeta} \varepsilon}{\chi},
\end{align}
is $\mathcal{O}(\log(\lceil 1/\epsilon \rceil))$.

For the case when errors constitute a white noise process, the following corollary provides the required result. 
\begin{cor}\label{lem:tracking_iid}
Under (\textbf{A1})-(\textbf{A5}), and sequences $\{\x_k^\star, \e_k\}$ satisfying \eqref{boundedekwk}, the IOGD iterates for $\pa$ adhere to the following bound
\begin{align}\label{tracking_per2}
\Ex{\dist{\x_{k+1},\cX_{k+1}^\star}}\leq &\left({{\ell}}/{\chi}\right)^k\dist{\x_1,\cX_1^\star}  \nonumber
\\
&+ \left[{\frac{1-\left({{\ell}}/{\chi}\right)^k}{1-\left({{\ell}}/{\chi}\right)}}\right]\frac{\sigma +\alpha \varepsilon}{\chi}.
\end{align}	
for all $k \geq 1$ and ${\ell}$ is as defined in Corollary \ref{p1iidcor} and $\alpha$ should be chosen in same manner as in Corollary \ref{p1iidcor}. 
\end{cor}
The proof of Corollary \ref{lem:tracking_iid}  \colb{follows from the analysis provided in Appendix B}.  Likewise, the iteration complexity of the tracking error for this case is also $\mathcal{O}(\log(\lceil 1/\epsilon \rceil))$.

\section{Implications of proposed results}\label{special}
This section provides two online algorithms that can be interpreted as variants of the proposed IOGD algorithm. 

\subsubsection{Incremental OGD with increasing sample size}\label{incre_IOGD}
Consider a scenario where the target trajectory \colb{can be expressed} as the solution to the following composite function optimization problem\vspace{-0mm}
\begin{align}
\x_k^\star &= \arg \min_{\x \in \cX} f_k(\x_k):= \frac{1}{N}\sum_{i=1}^N f_k^i(\x) & k = 1, 2, \ldots
\end{align}
where the objective function adheres to either $\pa$ or $\pb$. Inspired by the sampling-based gradient methods \cite{IGM_hy}, consider the following online algorithm:
\begin{align}
\x_{k+1} = \mathcal{P}_{\cX}\left(\x_k - \frac{\alpha}{N_k}\sum_{i\in \mathcal{N}_k} \nabla f_k^i(\x_k)\right) \label{iogdiss}
\end{align}
where $N_k = \abs{\mathcal{N}_k}$ and $\mathcal{N}_k$ is a random subset of $\{1, \ldots, N\}$. In order to view \eqref{iogdiss} as an IOGD variant, observe that the gradient error is given by $\e_k := \frac{1}{N_k}\sum_{i\in \mathcal{N}_k} \nabla f_k^i(\x_k) - \nabla f_k(\x_k)$. If the subsets $\mathcal{N}_k$ are formed by sampling the functions $f_k^i$ uniformly without replacement and in an i.i.d. fashion, it holds that [Sec. 2.8]{\cite{lohrsampling}},
\begin{align}
\Ex{\norm{\e_k}^2 \mid \F_k} = \frac{N-N_k}{N N_k} \Lambda^2
\end{align}
where $\Lambda^2$ is a bound on the sample variance of the gradients $\{\nabla f_k^i(\x_k)\}_{i=1}^N$, i.e., 
\begin{align}
\frac{1}{N-1}\sum_{i=1}^N \norm{\nabla f_k^i(\x) - \nabla f_k(\x)} &\leq \Lambda^2 & \x &\in \cX
\end{align}
Defining $\varepsilon_k:=\sqrt{(1/N_k-1/N)}$, it can be seen that $E_K$ is sublinear, for instance when $N_k=Nk^{\rho}/(N+k^{\rho})$ with $\rho >0$. 

\subsubsection{Proximal OGD methods for composite minimization}
Next, consider a scenario where the target trajectory is given by the solution of the following composite function minimization problem: 
\begin{align}
\x_t^\star = \arg\min_{\x} f_k(\x):=g_k(\x) + h_k(\x) \label{compos}
\end{align}
where $h_k$ is a differentiable regularization function. We propose the proximal OGD algorithm for solving \eqref{compos} that takes the form 
\begin{align}
\x_{k+1} = \text{prox}_{\alpha h_k}\left(\x_k - \alpha \nabla g_k(\x_k) \right) \label{prox}
\end{align}
where the proximal function is defined as $\text{prox}_{\alpha h_k}(\z) = \arg \min_{\x} \alpha h_k(\x) + \frac{1}{2}\norm{\x-\z}^2$. Equivalently, it is possible to write the update in \eqref{prox} as
\begin{align}
\nabla_{\x} \left(\alpha h_k(\x) + \frac{1}{2}\norm{\x - \x_k + \alpha \nabla g_k(\x_k)}^2\right) = 0 \nonumber\\
\Rightarrow \x_{k+1} = \x_k - \alpha \nabla g_k(\x_k) - \alpha \nabla h_k(\x_{k+1}) \label{pogd}
\end{align}
Recall that the IOGD algorithm for \eqref{compos} takes the form $\x_{k+1} = \x_k - \alpha (\nabla g_k(\x) + \nabla f(\x) +\e_k)$. Therefore, the error term becomes:
\begin{align}
\norm{\e_k} &= \norm{\nabla h_k(\x_k) - \nabla h_k(\x_{k+1})} \leq L \norm{\x_{k+1}-\x_k} \nonumber\\
&\leq \alpha L \norm{\nabla g_k(\x_k) + \nabla h_k(\x) + \e_k} \nonumber\\
&\leq \frac{\alpha L}{1-\alpha L} \norm{\nabla f_k(\x_k)}
\end{align} 
where $L$ is the Lipschitz constant of $\nabla h_k$. Consequently, as long as $\nu =\frac{\alpha L}{1-\alpha L}$ is sufficiently small and subject to the assumptions (\textbf{A1})-(\textbf{A5}), the proximal OGD algorithm adheres to the $\mathcal{O}(1+W_K)$ bound developed earlier. In other words, the inexact gradient used in \eqref{pogd} does not have any effect on the dynamic regret rate.

\vspace{-0mm}
\section{ Numerical Tests}\label{num}
This section provides detailed numerical tests that demonstrate the usefulness and applicability of the proposed IOGD framework to {three} different applications. 
\vspace{0mm}
%
%
%
		\subsection{Time varying parameter estimation} \label{time_varying_estimation}
			\begin{figure*}
				\setcounter{subfigure}{0}
				\begin{subfigure}{0.50\columnwidth}
					\includegraphics[width=\linewidth, height = 0.6\linewidth]
					{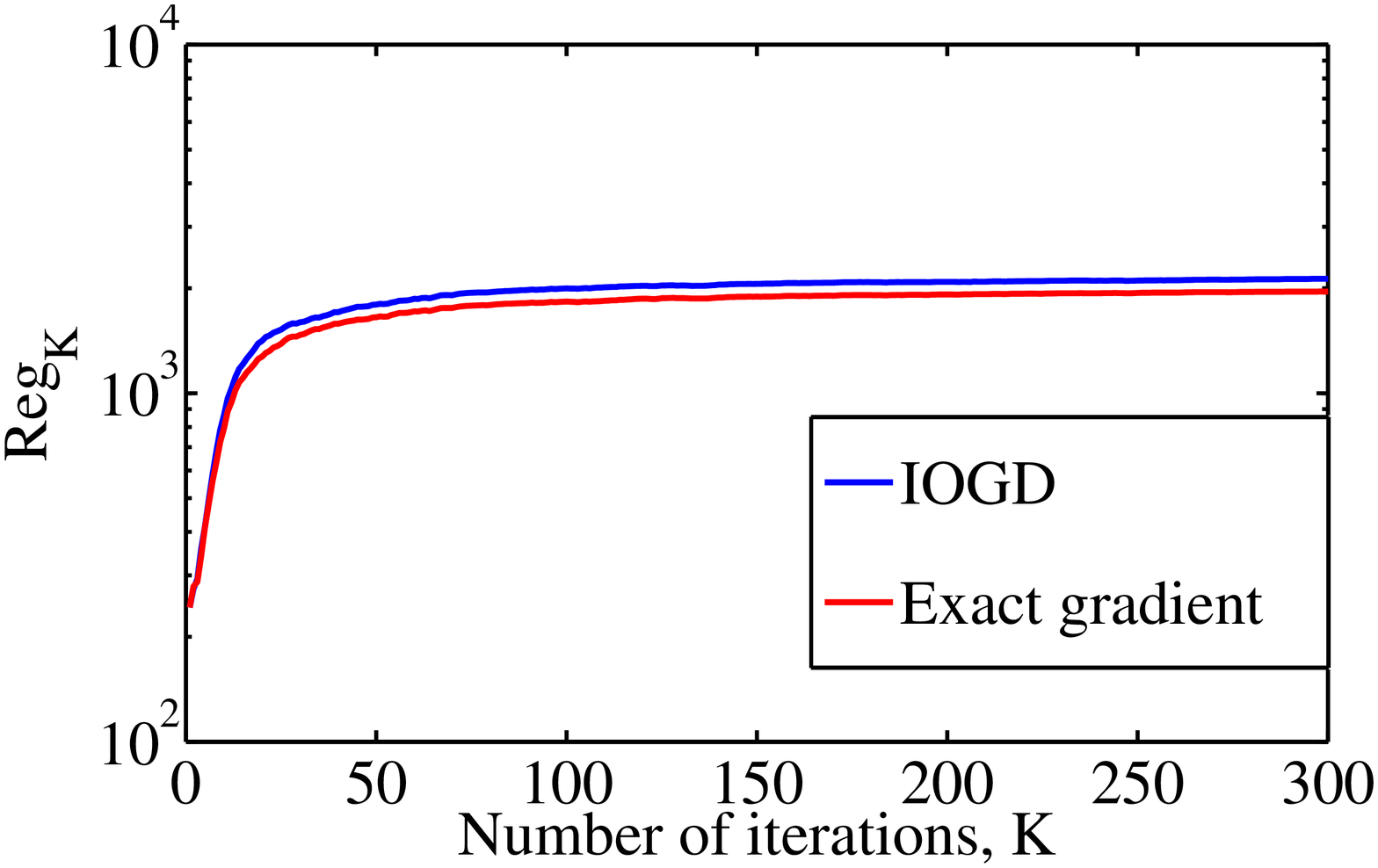}
					\caption{Dynamic regret performance}
					\label{est_1}
				\end{subfigure}
				\begin{subfigure}{0.50\columnwidth}
					\includegraphics[width=\linewidth,height = 0.6\linewidth]
					{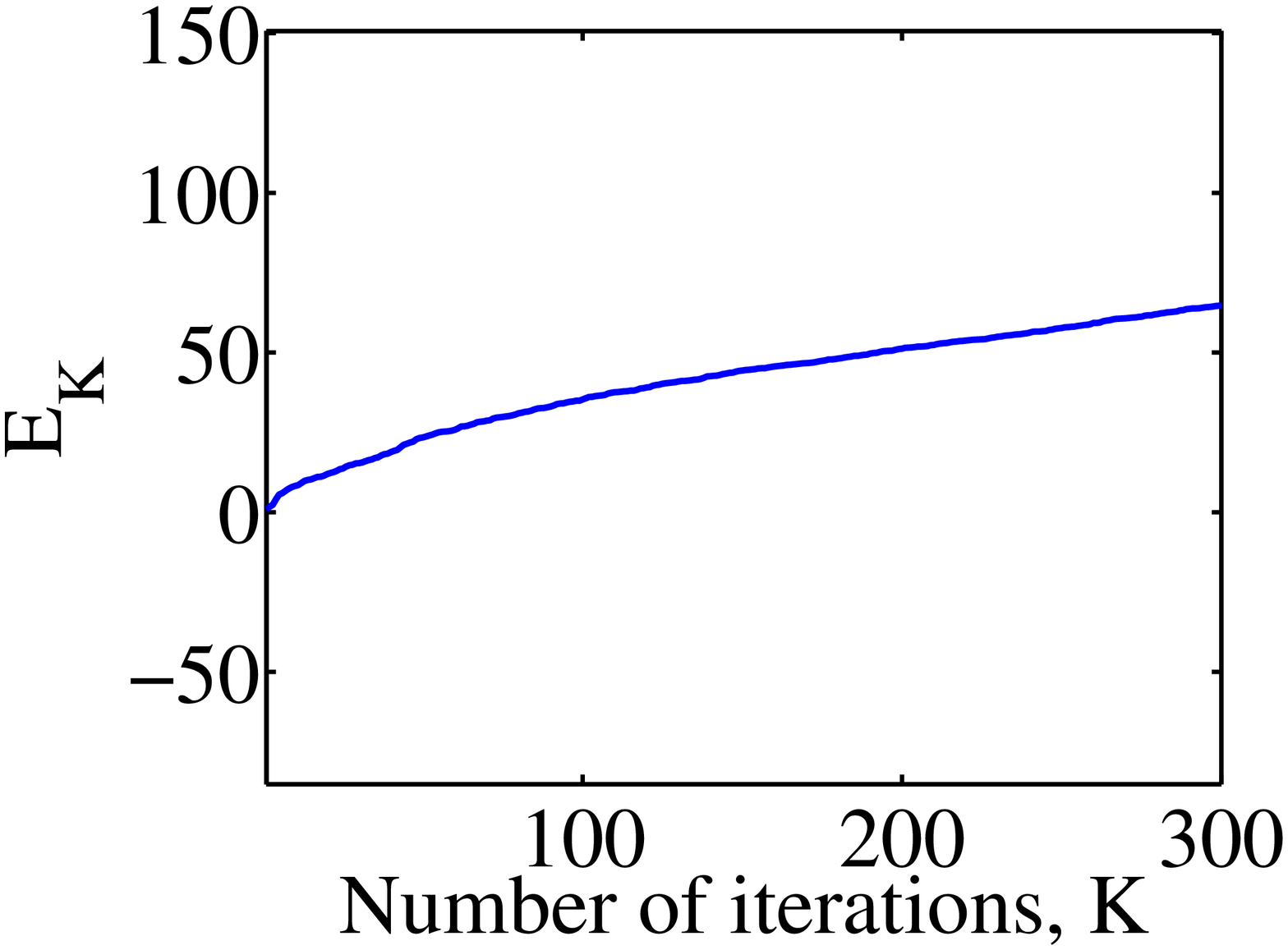}
					\caption{Cumulative error}
					\label{est_2}
				\end{subfigure}
					\begin{subfigure}{0.5\columnwidth}
						\centering
						\includegraphics[width=\linewidth,height = 0.6\linewidth]
						{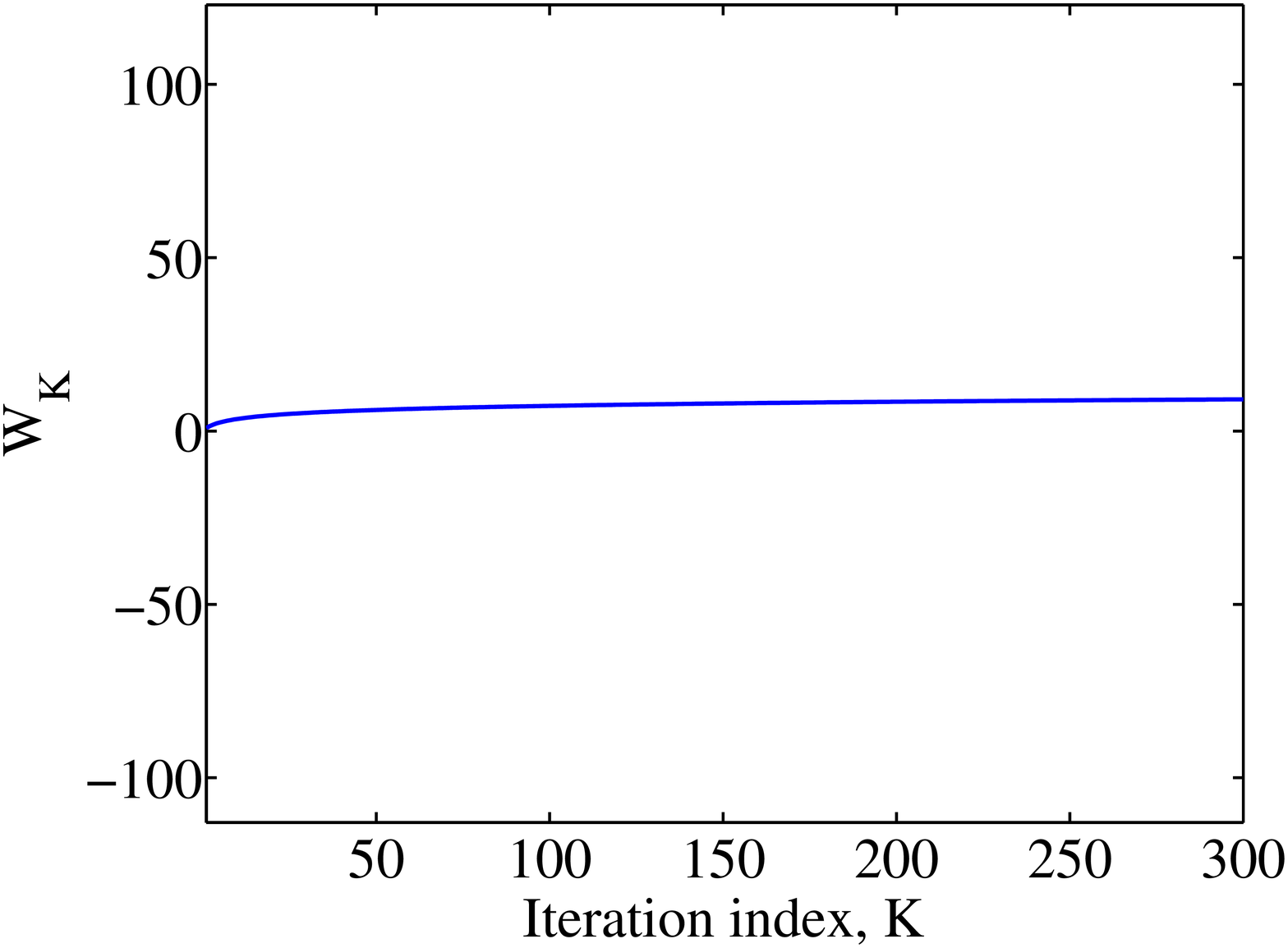}
						\caption{Path length}
						\label{est_3}
					\end{subfigure}
				\caption{\colb{Time varying parameter estimation for $m=2$ and $n=3$, $N=100$. Fig. 1(a) plots the dynamic regret, Fig. 1(b) plots the cumulative error $E_K$, and Fig. 1(c) shows the sublinear behavior of path length $W_K$ as a function of the total number of iterations $K$.}}\vspace{0mm}
			\end{figure*}
	\colb{	Consider a network of $N$ nodes placed over a planar region $\mathcal{A}\subseteq \mathbb{R}^2$. The nodes seek to cooperatively estimate (track) a time-varying parameter $\x^\star_k\in\Rn^n$ using observations $\b_{i,k}\in\Rn^m$ made at each node $i$ and time $k$. The measurements follow the linear model $\b_{i,k}=\E_{i}\x_{k}^\star+ \n_{i,k}$ where $\E_i$ is a node-dependent measurement matrix and the noise $\n_{i,k} \sim \mathcal{N}(0, \sigma^2 \boldsymbol{I})$ is independent identically distributed across nodes and time. At time $k$, the least squares estimate for $\x_k$ is given by 
\begin{align} \label{eq:dist_rls}
\x_k^\star=\arg\min_{\x\in \Rn^n}\frac{1}{N}\sum_{i=1}^{N} \norm{\E_{i}\x_{k} - \b_{i,k}}^2_2. 
\end{align} 
Here, the measurement matrices $\{\E_i\}_{i=1}^N$ may not necessarily be full rank and therefore the optimum $\x_k^\star$ may not be unique. In this case, the path length in \eqref{pl} is well-defined for a given trajectory $\{\x_k^\star\}$, while the path lengths in \eqref{rem_first} and \eqref{rem_second} turn out to be infinite.} 

\colb{Towards estimating $\x_k$ in an efficient manner, we make use of the incremental OGD algorithm proposed in \eqref{iogdiss}. Letting $\mathcal{N}_k$ be a random subset of $\{1, \ldots, N\}$, the incremental OGD updates take the form:
\begin{align}
\x_{k+1} = \x_k - \frac{\alpha}{N_k}\sum_{i\in \mathcal{N}_k} \nabla f_k^i(\x_k) \label{numerical_1}
\end{align} 
where $f_k^i(\x_k) := \norm{\E_{i}\x_{k} - \b_{i,k}}^2_2$ and $N_k := \abs{\mathcal{N}_k}$. As detailed in Sec. \ref{special}, the incremental OGD is a special case of the proposed IOGD algorithm if $N_k$ increases to $N$ as $k$ goes from 1 to $K$. Observe that the gradient  error for this case is not independent or identically distributed. Therefore, for this case, assumptions (\textbf{A1}) and (\textbf{A3})-(\textbf{A5}) are satisfied and the results of Sec. \ref{rbp1} are applicable.}

\colb{In order to validate the performance of the proposed algorithm, we consider a network with $N = 100$ nodes and for $m=2$ and $n=3$. The subset $\mathcal{N}_k$ is selected uniformly without replacement from $\{1, \ldots, N\}$ and following the rule $N_k = Nk/(N+k)$ as stated in Sec. \ref{special}. The step size is set to $\alpha = 0.1$ for all iterations and irrespective of $K$. The algorithm is run for different values of $K$ and the resulting dynamic regret, error bound, and the path length are shown in Figs. 1(a), 1(b), and 1(c). For the purposes of comparison, the regret obtained when full data is used (i.e. $N_k = N$ for all $k\geq 1$) is also shown (exact gradient). It can be seen that although proposed incremental IOGD does not make use of the full data, its regret is not too far from that of the classical OGD algorithm. Further, as expected from the analytical results, the dynamic regret is clearly sublinear. }


\subsection{Multi-target tracking}\label{SAP} 
This section develops a low-complexity online multi-target tracking algorithm inspired from the convex optimization-based target tracking framework developed in \cite{derenick2009convex}. Specifically, a team of $n$ agents at locations $\{\x^i_k\}_{i=1}^n$ is tasked with tracking a set of $m$ targets at locations $\{\y_k^j\}_{j=1}^m$. The discretized problem is formulated as the following convex optimization that must be solved for each $k\geq 1$ \cite[Thm. 3.8.1]{derenick2009convex}:
\begin{subequations}\label{socp}
\begin{align}
\{\x_{k+1}^i\} = &\arg\min_{\{\x^i\in\Rn^p\}}  \sum_{i=1}^n \psi_k^i(\x^i) &\\
		&\text{s.t.}\   \norm{\x^i-\x_k^i}^2 \leq v \ \ ;\ \   i=1,\ldots,n\label{cons}\\
		& \sum_{i=1}^n w^{ij}_k \norm{\x^i-\y_k^j}^2 \leq \eta	\ \ ;\ \  j = 1, \ldots, m \label{wcon}
		\end{align}
\end{subequations}
where $\psi_k^i(\cdot)$ is a time-varying cost function and $v$ is the square of the maximum distance that an agent may cover within a single time slot. For this paper, the following cost function is used:
\begin{align}
\psi_k^i(\{\x^i\})= \frac{1}{2}\sum_{\ell\neq i} \norm{\x^i -\x^\ell}^2 + \gamma  \norm{\x^i -\x^i_k}^2 \label{cost}
\end{align}
	where $\gamma >0$ is a regularization parameter. The objective function encourages agent $i$ to remain close to the other agents. At the same time, the regularization term forces the agents to not move around unnecessarily. The constraint in \eqref{wcon} is the linearized version of the original constraint obtained from using the process described in \cite[Chap. 3]{derenick2009convex}. A sigmoidal weight function is utilized that takes the form:
\begin{align}
w_k^{ij} = \left(1+e^{-\vartheta(\epsilon-\norm{\x^i_k-\y_k^j})}\right)^{-1} \label{weight}
\end{align}
where $\vartheta$ and $\epsilon$ are positive parameters. Observe that the weights are small for agent-target pairs that are far from each other. In other words, the constraint in \eqref{wcon} encourages the set of agents tracking a target $j$ to stay close to it. The weights are normalized so that {$\sum_{i=1}^nw^{ij}_k = 1$} in order to ensure that each target is tracked by at least one agent. Finally, it is remarked that the agents may only know the estimated target locations $\{\hat{\y}_k^j\}_{j=1}^m$ instead of the true locations required for \eqref{wcon}. 

Different from \cite{derenick2009optimal}, where interior point methods are proposed towards solving \eqref{socp} at each time instant, the goal here is to track the targets at a significantly low complexity. To this end, we propose an IOGD-based low-complexity algorithm capable of tracking high speed targets in large multi-agent networks. Since the constrained optimization problem in \eqref{socp} is not of the form required in \eqref{xkstar}, the IOGD algorithm will instead be applied in the dual domain. It can be verified that \eqref{socp} satisfies the Slater's conditions, and therefore has zero duality gap.


To this end, associate dual variables $\{\lambda^i\}_{i=1}^n$ and $\{\nu^j\}_{j=1}^m$ with \eqref{cons} and \eqref{wcon}, respectively. Collecting the primal and dual variables $\{\lambda^i\}$, $\{\nu^j\}$, and $\{\x^i\}$ into vectors $\lb$, $\nb$, and $\x$ respectively, the Lagrangian can be written as
\begin{align}
L_k(\x,\lb,\nb) &= \sum_{i=1}^n\sum_{\ell = i+1}^n \norm{\x^i-\x^\ell}^2 + \gamma \norm{\x-\x_k}^2 \\
&\hspace{-1.8cm} + \sum_{i=1}^n \lambda^i(\norm{\x^i-\x_k^i}^2-v) + \sum_{j=1}^m\nu^j\left[\sum_{i=1}^n w_k^{ij}\norm{\x^i-\y^j_k}^2 - \eta\right].\nonumber
\end{align}
Thus the dual function can be written as
\begin{align}
\varrho_k(\lb,\nb) = \arg \min_{\x} L_k(\x,\lb,\nb)
\end{align}
Since the dual function is always concave, the proposed IOGD algorithm can be utilized to maximize  $\varrho_k$ in an online fashion. In particular, for a given $\lb$ and $\nb$, the gradient of the dual function is given by 
\begin{align}
\frac{\partial L_k(\x,\lb,\nb)}{\partial \lambda^i} &= \norm{\x^i(\lb,\nb)-\x^i_k}^2-v \label{gradlam}\\
\frac{\partial L_k(\x,\lb,\nb)}{\partial \nu^j} &= \sum_{i=1}^n w_k^{ij}\norm{\x^i(\lb,\nb)-\y^j_k}^2-\eta \label{gradnu}
\end{align}
where $\{\x^i(\lb,\nb)\}_{i=1}^n = \arg \min_{\x} L(\x,\lb,\nb)$. In order to evaluate the minimization, observe that the condition $\nabla_{\x^i} L_k(\x,\lb,\nb) = 0$ is equivalent to
\begin{align}
&(n + 2\gamma + 2\lambda^i + 2\sum_{j=1}^m \nu^jw_k^{ij})\x^i - \sum_{\ell=1}^n \x^\ell \nonumber\\
&\hspace{3cm} = 2(\gamma+ \lambda^i)\x_k^i + 2\sum_{j=1}^m \nu^jw_k^{ij}\y_k^j \label{Lgrad}
\end{align}
for $i =1, \ldots, n$. Let $\D_k$ be a diagonal matrix with $(i,i)$-th element $[\D_k]_{ii}=n + 2\gamma + 2\lambda^i + 2\sum_{j=1}^m \nu^jw_k^{ij}$ and let the right-hand side of \eqref{Lgrad} be denoted by the $i$-th element of the vector $\boldsymbol{\varphi}_k$. Then the solution to \eqref{Lgrad} can be written as $\x = (\D_k-\mathbf{11}^T)^{-1}\boldsymbol{\varphi}_k$ where $\mathbf{1}$ is an all-one vector of appropriate dimension. Observe further that the solution is unique for all $\lb$ and $\nb$ since $L_k(\x,\lb,\nb)$ is strongly convex in $\x$. Using the matrix inversion lemma, it is possible to calculate $\x$  with $\mathcal{O}(n)$ complexity as follows:
\begin{align}\label{primal}
\!\!\!\!\x^i\!\! =\!\! [\D_k]^{-1}_{ii}[\boldsymbol{\varphi}_k]_i\!\! +\! \frac{[\D_k]^{-1}_{ii}}{1\!\!-\!\sum_{\ell=1}^n [\D_k]_{\ell\ell}^{-1}} \left(\sum_{\ell=1}^n [\D_k]^{-1}_{\ell\ell}[\boldsymbol{\varphi}_k]_\ell\right)
\end{align}
Since  the target location is not known exactly, the estimated target location must be used in \eqref{gradnu}, making the gradient inexact. \colb{The agent location is determined from the current primal iterate $\x_k$ which serves as a proxy for the primal solution $\x_k^\star$ at time $k$. Consequently, the results of Sec. \ref{track} will be applicable here. } The full IOGD-based multi-target tracking algorithm is summarized in Algorithm \ref{Algo_1}. 
 
\begin{figure}[h]
 		\centering
 		\includegraphics[width=\linewidth, height = 0.6\linewidth]{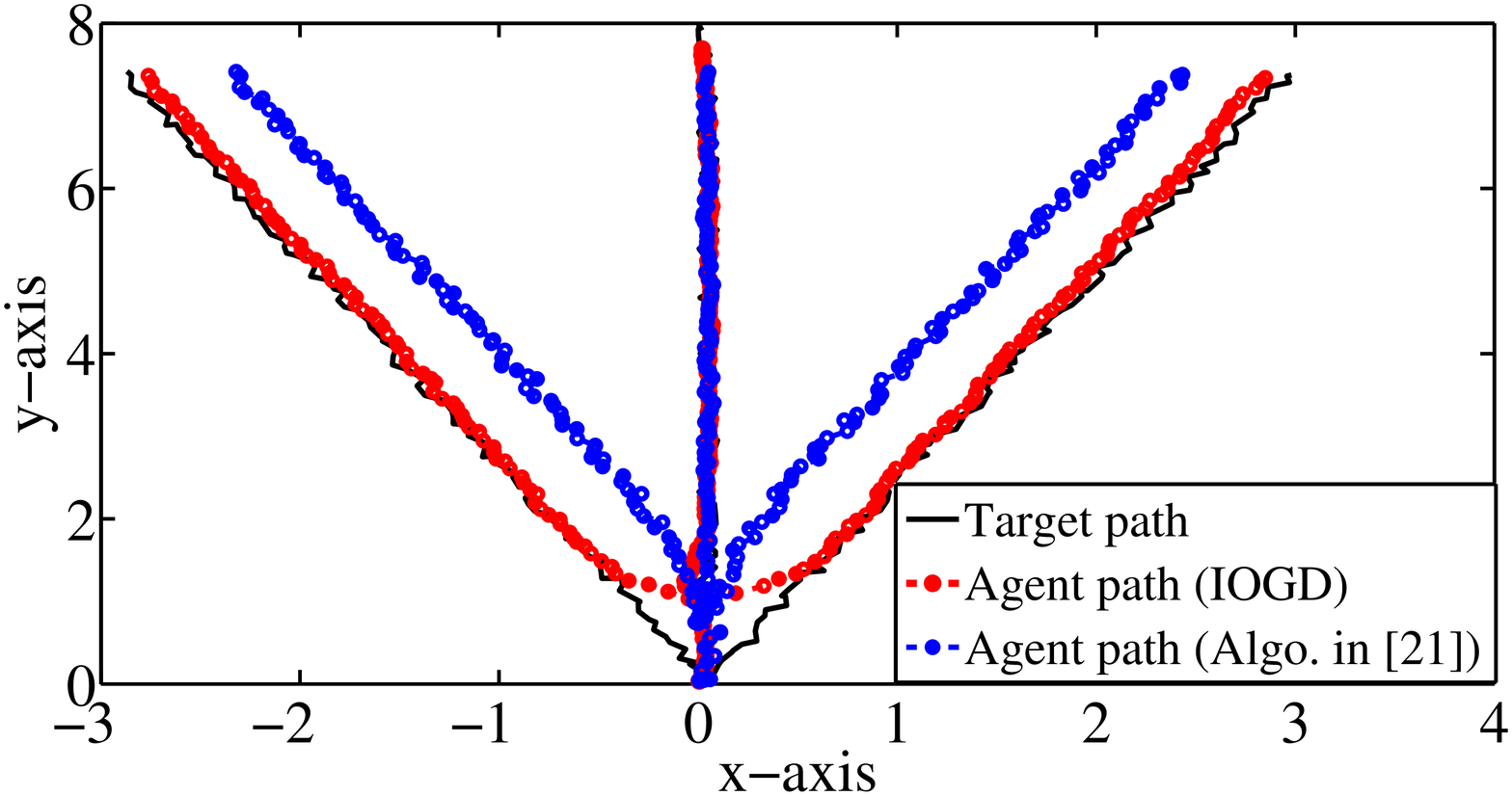}
 		\caption{Tracking performance comparison for $m=3$ and  $n=3$.}
 		\label{Target_tracking}
 	\end{figure}

\begin{algorithm} 
		\caption{IOGD-based multi-target tracking}
		\begin{algorithmic}[1]\label{Algo_1}
			\STATE {\bf Initialize} $\x_1$, $\lb_1$, and $\nb_1$, and step sizes $\alpha_{\lambda}$, and $\alpha_{\nu}$  \\
			\textbf{Repeat for $k = 1, 2, \ldots, $}
			\STATE \hspace{0.5cm} Compute weights $\{w_k^{ij}\}$ from \eqref{weight}\\
			\STATE \hspace{0.5cm} \textbf{Calculate} the next location as \\ 
			\hspace{0.5cm} $\x_{k+1} = \arg\min_{\x} L_k(\x, \lb_k, \nb_k)$ 
			\STATE \hspace{0.5cm} Update for all agents and targets:
			\begin{align}
			\lambda^i_{k+1} &= \left[\lambda^i_{k} - \alpha_{\lambda}(\norm{\x^i_{k+1}-\x^i_k}^2-v)\right]_{+} \label{iglam}\\
			\nu^j_{k+1} &= \left[\nu^j_{k} - \alpha_{\nu}(\sum_{i=1}^n  w_k^{ij}\norm{\x^i_{k+1}-\hat{\y}^j_k}^2-\eta)\right]_{+} \label{ignu}
			\end{align}
		\end{algorithmic}
	\vspace{0mm}
	\end{algorithm} 	
	\vspace{0mm}
The performance of the proposed multi-target tracking algorithm is studied on a number of simulated planar environments. \colb{The agent velocities are restricted to 0.2 m/s and a target is assumed covered if it is within $\eta=0.54$ m from at least one agent.} Consider first a simple scenario consisting of three targets ($m=3$) and three agents $n=3$. The targets are co-located at time $k=1$ and start moving away from each other along the paths shown in Fig. \ref{Target_tracking}. It can be seen that the proposed algorithm works as expected, and the agent team splits up in order to track the three targets. On the other hand, the algorithm in \cite{derenick2009optimal} does not necessarily exhibit such a behavior and requires careful parameter tuning so as to allow tracking with reasonable accuracy; see Fig. \ref{Target_tracking}. Indeed, since \cite{derenick2009optimal} entailed  solving a constrained convex optimization problem at every time instant, it was observed that unless the parameters are not selected carefully, the problem could become infeasible. It was however possible to circumvent this behavior to a certain extent by explicitly adding noise to the output of the optimization problem. In contrast, no such issue was present in the proposed IOGD algorithm, whose performance was quite robust to the choice of parameters. \colb{It is observed from the figure that the trajectory of the solution of (50) using [21] is a bit away from the target trajectory because the formulation in (50) does not incentivize the agents to come too close to target. On the other hand, the IOGD algorithm over-tracks to a certain extent, thereby coming close to the target but at the same time yielding a higher objective value.}

\begin{figure}[h]
		\centering
		\includegraphics[width=\linewidth, height = 0.4\linewidth]{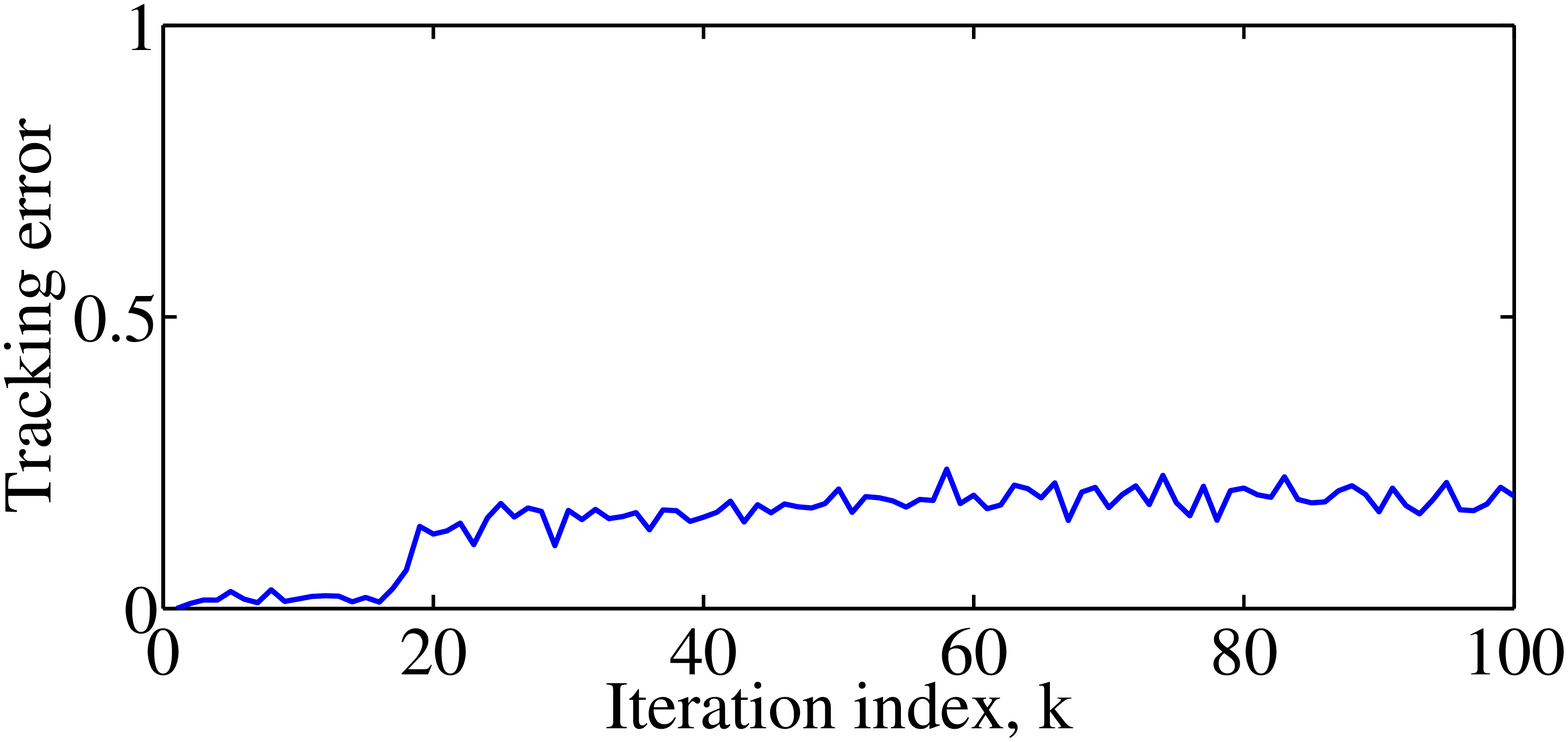}
		\caption{Tracking error performance of IOGD algorithm. }
		\label{error_bound}
	\end{figure}

\colb{ Fig. \ref{error_bound} shows the error performance $\norm{\x_k-\x_k^\star}$ for the IOGD algorithm. As expected from the results of Sec. \ref{track}, the tracking error remains bounded. It is remarked that for tracking error is small initially while the agents are close to each other and increases when they split up to track the diverging targets.}
 
%
	
\begin{figure*}
		\setcounter{subfigure}{0}
		\begin{subfigure}{0.5\columnwidth}
			\includegraphics[width=\linewidth, height = 0.6\linewidth]
			{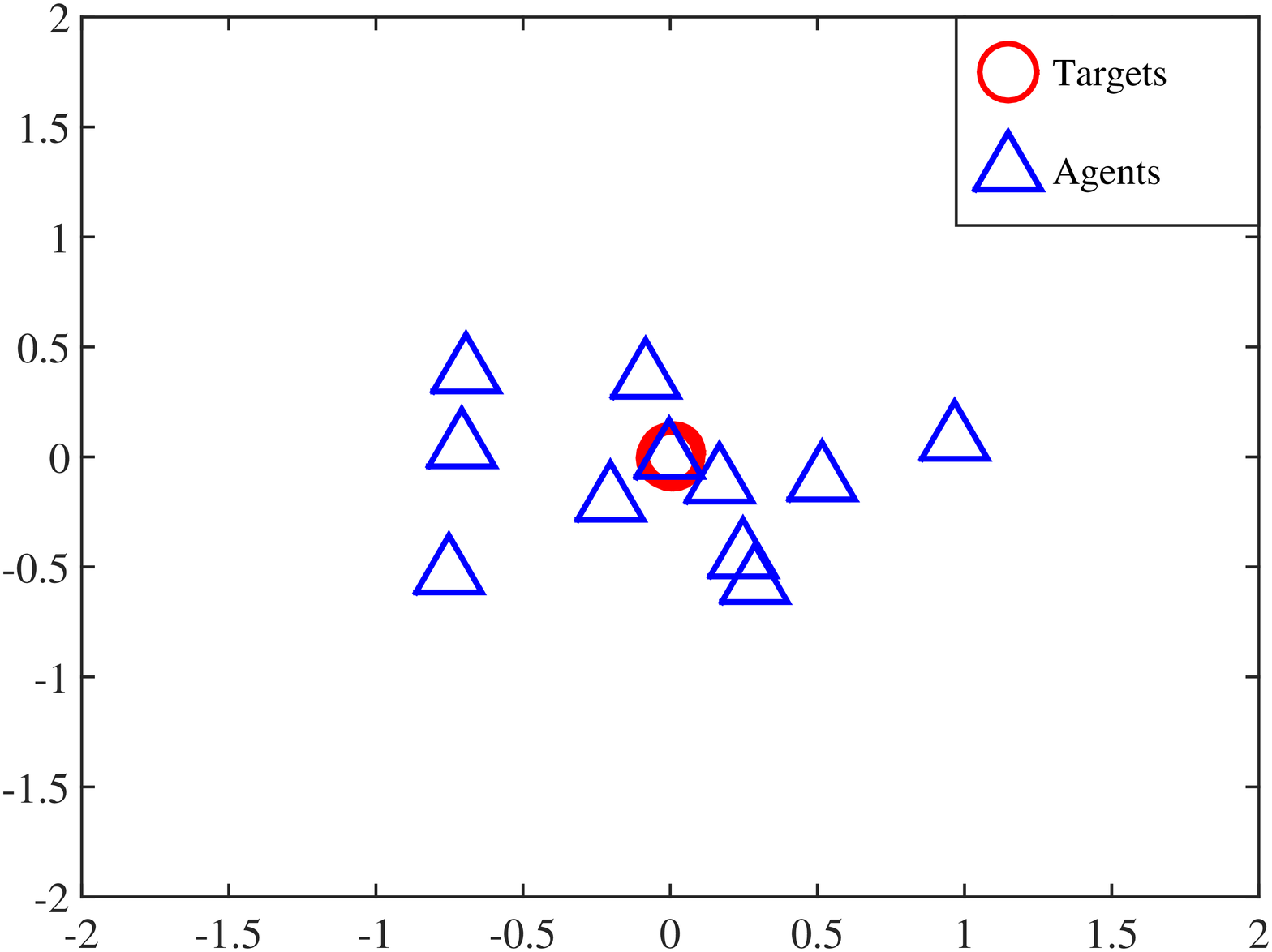}
			\caption{$k=1$}
			\label{fig0}
		\end{subfigure}
		\begin{subfigure}{0.5\columnwidth}
			\includegraphics[width=\linewidth,height = 0.6\linewidth]
			{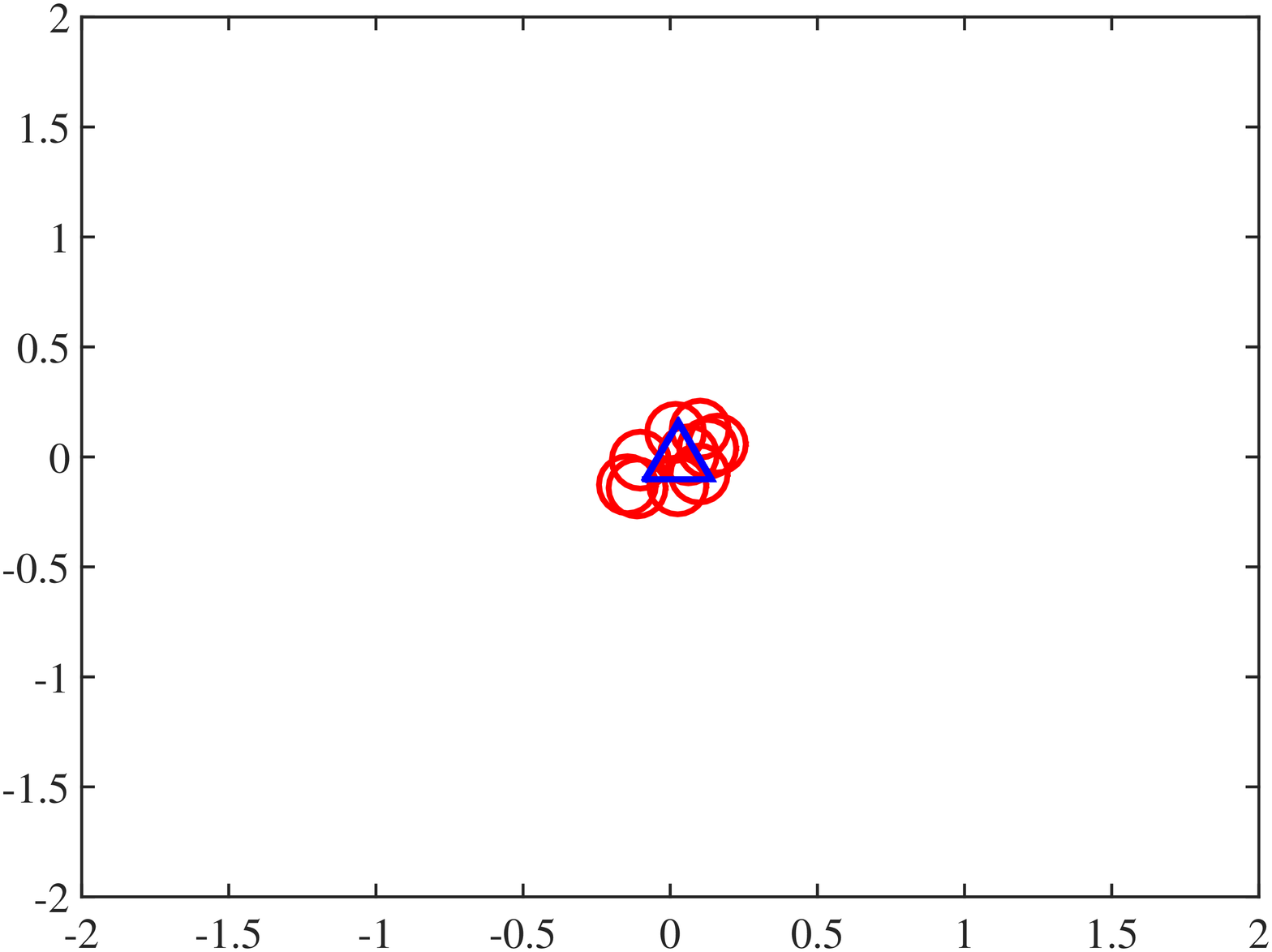}
			\caption{$k=10$}
			\label{fig1}
		\end{subfigure}
		\begin{subfigure}{0.5\columnwidth}
			\includegraphics[width=\linewidth,height = 0.6\linewidth]
			{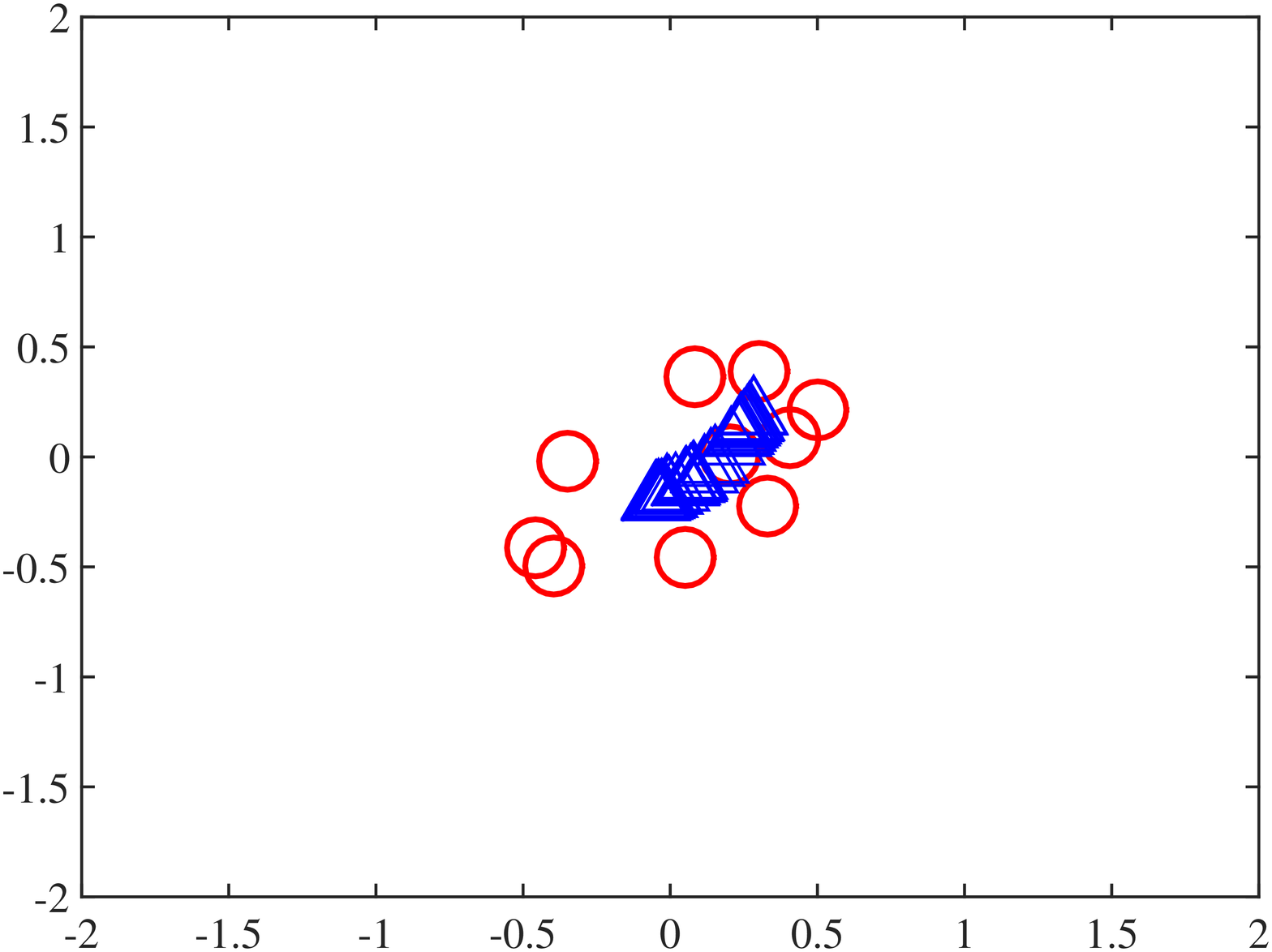}
			\caption{$k=30$}
			\label{fig2}
		\end{subfigure}
		\begin{subfigure}{0.5\columnwidth}
			\includegraphics[width=\linewidth, height = 0.6\linewidth]
			{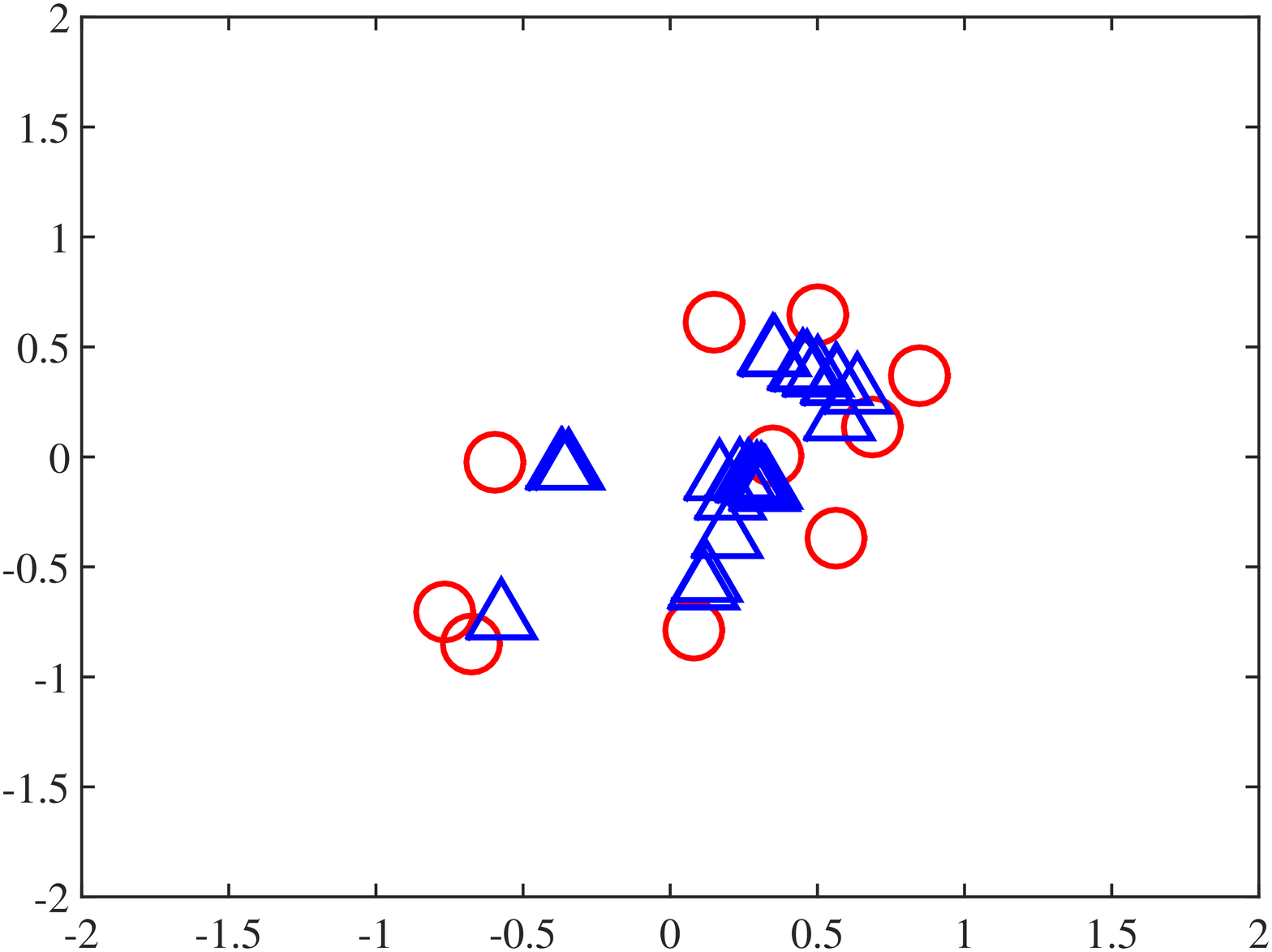}
			\caption{$k=50$}
			\label{fig3}
		\end{subfigure}
		\begin{subfigure}{0.5\columnwidth}
			\includegraphics[width=\linewidth,height = 0.6\linewidth]
			{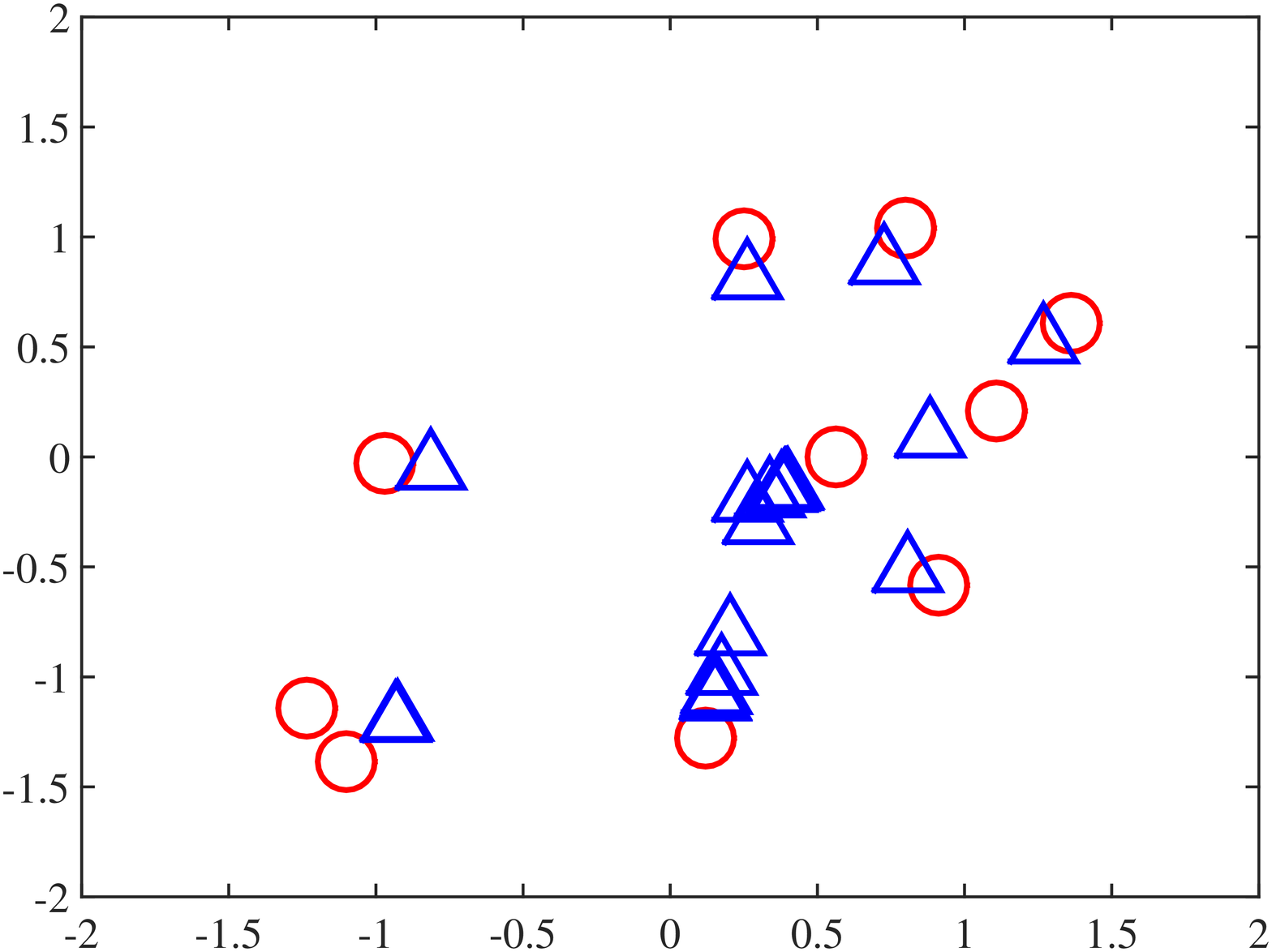}
			\caption{$k=80$}
			\label{fig4}
		\end{subfigure}
		\begin{subfigure}{0.5\columnwidth}
			\includegraphics[width=\linewidth,height = 0.6\linewidth]
			{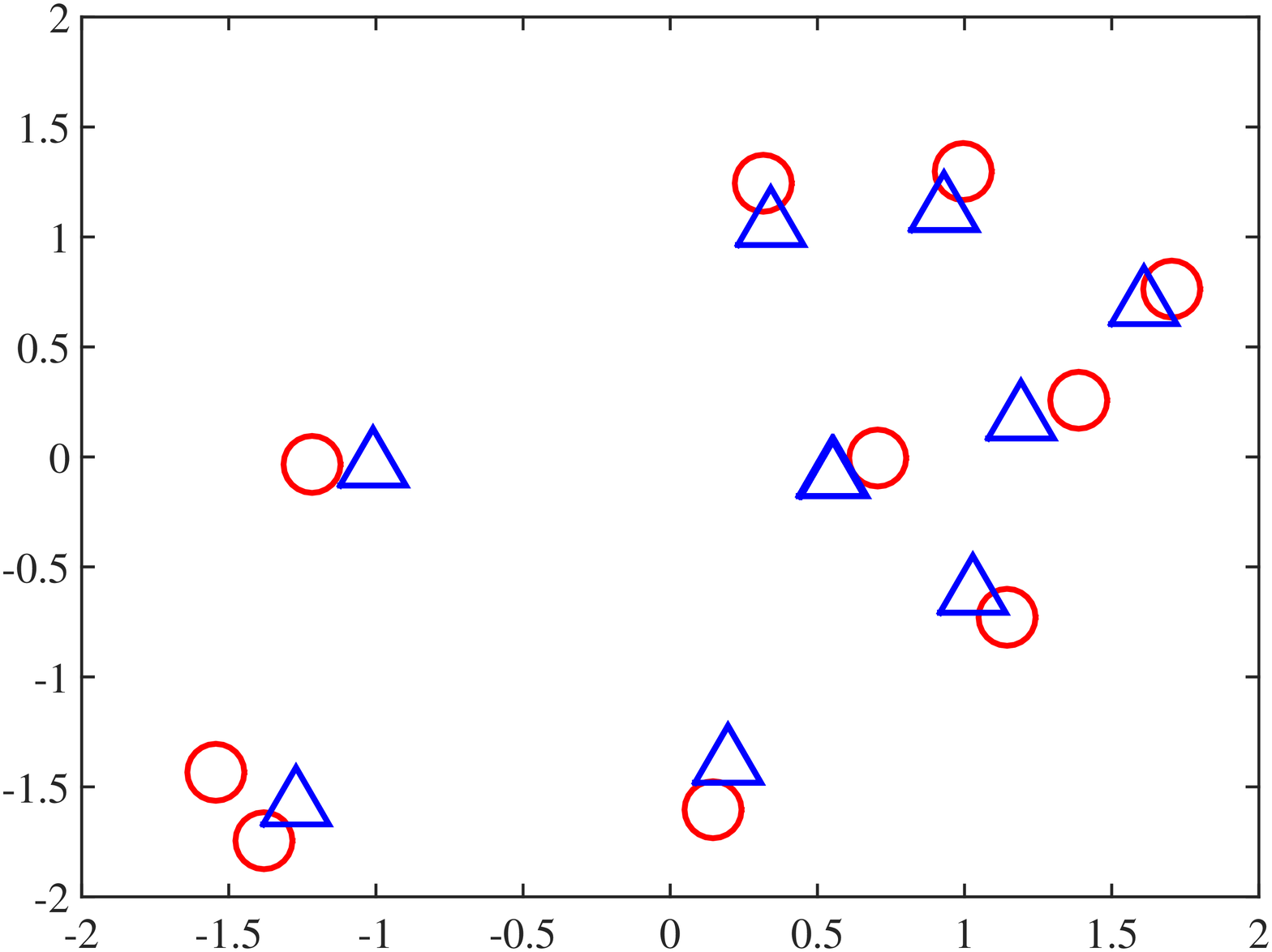}
			\caption{$k=100$}
			\label{fig5}
		\end{subfigure}
		\caption{Target tracking for $m=10$ and $n=50$.}\label{atc2}\vspace{0mm}
	\end{figure*}

\begin{figure}
		\centering
		\includegraphics[width=\columnwidth, height=0.5\textheight]{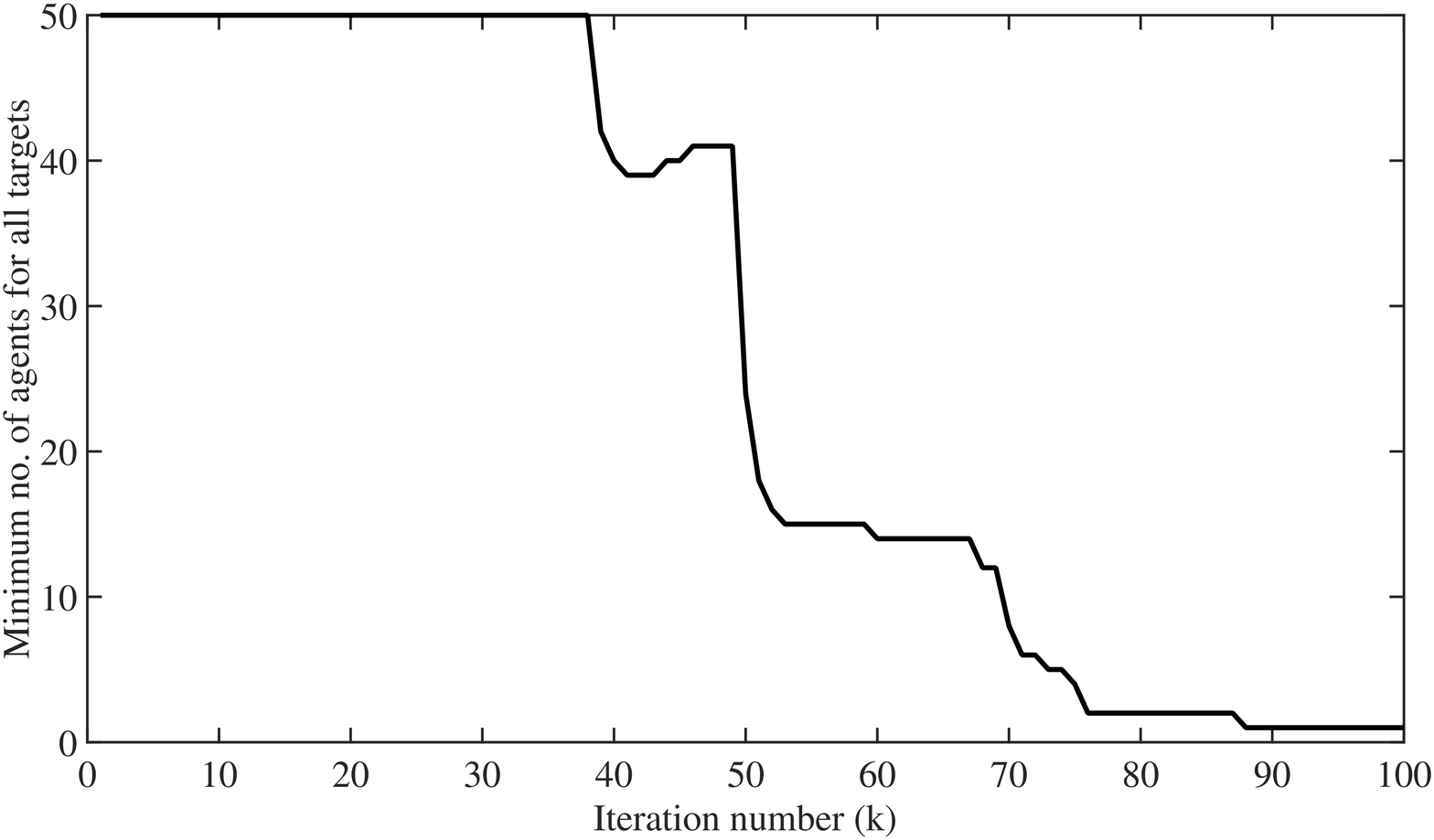}
		\caption{Number of agents covering each target ($m = 10$, $n = 50$).}
		\label{atc}
	\end{figure}

Next, we consider a large scale system with $m = 10$ targets and $n = 50$ agents. As expected, the IOGD algorithm is capable of tracking most of the targets at low complexity; see Figs. \ref{atc2} and \ref{atc}. As with the smaller system considered earlier, the splitting of the agent teams is observed in Fig. \ref{fig3}. The splitting behavior is also evident from the supplementary video included with this paper \footnote{\url{https://www.youtube.com/watch?v=bVto6LItehM}}. It is important to emphasize that the tracking performance of the IOGD is at par with the convex optimization approach of \cite{derenick2009optimal}. In contrast, solving a general convex optimization problem as required in \cite{derenick2009optimal} incurs a complexity of at least $\mathcal{O}(n^3)$ as opposed to the $\mathcal{O}(n)$ complexity incurred in the calculation of the inexact gradient in \eqref{iglam}-\eqref{ignu}. For the sake of comparison, both algorithms were implemented in MATLAB and their run-times measured on an Intel Xeon E3-1226 3.30GHz CPU machine. The resulting per-iteration run-time for the proposed algorithm was $49$ ms, as compared to that of $974$ ms required by \cite{derenick2009optimal}.

\subsection{ Tracking movie ratings }
This section considers the problem of tracking movie ratings in recommender systems. Traditionally, the problem has been posed as that of completing a low-rank matrix and solved using a variety of offline or incremental algorithms. In many applications however, the dataset grows over time, necessitating a real-time recommendation engine that predicts the missing ratings in an online manner; see e.g. \cite{dhanjal2014online}. This section advocates an IOGD-inspired algorithm that can predict and track missing ratings at low complexity. 

Let $\M_k \in \Rn^{K \times L}$ be the measurement matrix containing user ratings at locations specified in the binary matrix $\J_k$, i.e., $[\J_k]_{ij} = 1$ if $[\M_k]_{ij}$ contains a rating and zero otherwise. The goal is to generate a complete ratings matrix $\X_k$, which is the solution to the following time-varying optimization problem:
\begin{align}\label{movie_lens}
\X_k^\star = \min_{\X\in\cX} \frac{1}{2}\norm{\J_k\odot(\M_k-\X)}^2_F+\lambda\norm{\X}_{\star}.
\end{align}
Following the notation of \cite[Sec. III-D]{7903718}, $\odot$ denotes the Hadamard product while $\norm{\X}_\star$ denotes the nuclear norm of $\X$. In order to obtain $\X_k$ in an online fashion, we utilize the proximal OGD algorithm introduced in \eqref{prox}, that takes the form:
\begin{align}\label{movielens_algo}
\X_{k+1}=D_{\lambda}\left(\X_k+\mu\J_k\odot(\M_k-\X)\right)
\end{align}
where, $D_{\lambda}(\cdot)$ is the singular value thresholding operator. Given the singular value decomposition of a matrix $\Y = \U\text{diag}(\sigma_1, \ldots, \sigma_r)\V^T$, the singular value thresholding operation is defined as
\begin{align}
D_{\lambda}(\Y)=&\U \text{diag}([\sigma_1-\lambda]_{+}, \ldots, [\sigma_r-\lambda]_{+})\V^T.
\end{align}
It is remarked that although the proximal OGD in Sec. \ref{special} is motivated as an IOGD variant for differentiable regularizers, it is applied and will be shown to work even for the non-differentiable regularizer $\norm{\X}_\star$. \colb{The other assumptions, namely (\textbf{A1}), (\textbf{A4}), and (\textbf{A5}) are however satisfied for this case. } 

Towards testing the performance of the proposed algorithm, we consider the Movielens 10M dataset, consisting of about 10 million ratings. The time-stamp associated with each rating is utilized to obtain 40 matrices $\{\M_k\}_{k=1}^{40}$, each consisting of new ratings added over the last 30 days. Fig. \ref{mov_plot} shows the evolution of the root mean-square error over time for the step size parameter $\alpha = 1$. {Interestingly, the regularized mean square error (RMSE) decays quickly over time and attains a value competitive with the benchmark \cite{dhanjal2014online}. Indeed, the proximal algorithm in (59) can be viewed as a generalization of the Soft-Impute method for the dynamic matrix completion problem. The tracking MSE for a class of related adaptive matrix completion algorithms is provided in \cite{7903718}.

\section{Conclusion and Future Work}\label{Conclusion}
This paper considered the problem of tracking time-varying and possibly adversarial targets. An inexact gradient descent method is proposed, that is applicable to non-strongly convex objective functions. The performance of the proposed framework is analyzed by developing bounds on its dynamic regret in terms of the target path length and the cumulative gradient error. Two distinct cases are considered: (a) case $\pa$, where the cost function adheres to the quadratic growth condition; and (b) case $\pb$ where the cost function is convex but the optimization variable belongs to a compact domain. For both cases, it is established that the dynamic regret is sublinear if the target path length and the cumulative gradient errors are {sublinear}. Steady state tracking errors are developed for the scenarios where these quantities are not sublinear. Further, the proposed framework is utilized to develop online variants of the incremental gradient and the proximal gradient descent algorithms. Tests on the multi-target tracking and movie recommendation problems showcase the efficacy and low-complexity of the proposed algorithm. 

\colb{While the present work broadens the applicability of OGD-like methods to time-varying non-strongly convex functions, it also motivates some interesting research directions. To begin with, the assumptions made here still preclude a large class of machine learning problems with non-differentiable objective functions. More generally, it may be possible to relax (\textbf{A4}) and characterize the dynamic regret for general time-varying and constrained convex optimization problems. In practical multi-agent systems, it may also be desirable that the nodes minimize their communication overhead and use distributed OGD variants such as the cyclic incremental gradient descent method.  Next, the existing literature on online learning relies heavily on characterizing the performance of the algorithms using the dynamic regret as defined in \eqref{reg}. In the context of target tracking however, it may be useful to consider a more appropriate exponentially weighted dynamic regret, defined as $\reg^\omega:=\sum_{k=1}^{K}\omega^{K-k}(f_k(\x_k)-f_k(\x_k^\star))$ for $0 < \omega < 1$. By assigning higher weights to the decisions made in the recent past, such a definition is apposite for tracking problems with shorter memory. 
Such an exponential window function is also reminiscent of the classical adaptive filtering algorithms such as the recursive least squares algorithm. }

\begin{figure}
	\centering
	\includegraphics[scale=0.5]{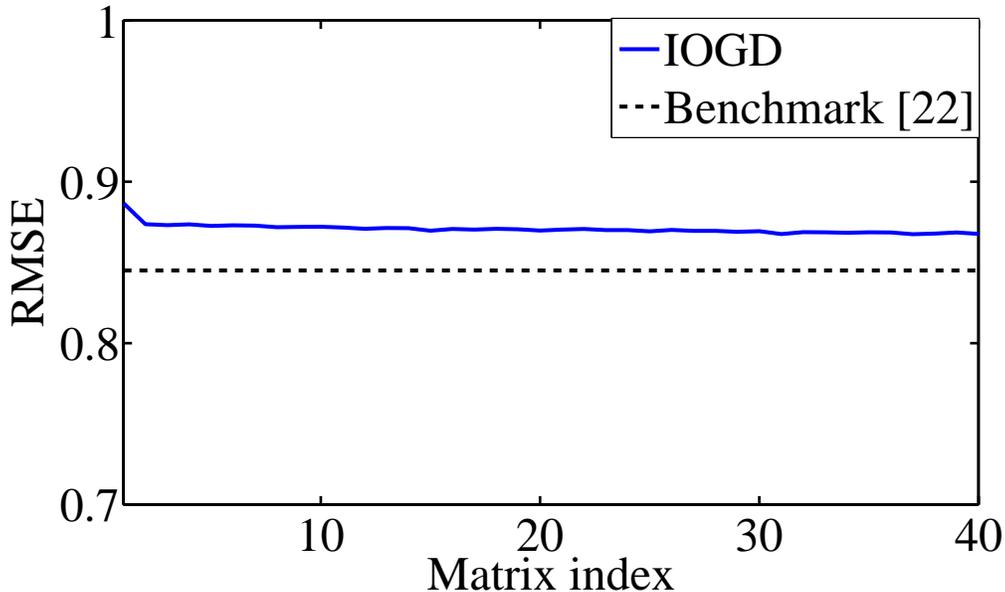}
	\caption{Tracking RMSE for the MovieLens Dataset}
	\label{mov_plot}
\end{figure}

	\appendices 
\section{Preliminaries}\label{prelim}
To begin with, several preliminary results required for establishing various bounds are derived. Observe that for $\pa$, it follows from the convexity and quadratic growth property of $f_k$ that for all $\x \in \cX$,
\begin{align}
\nabla f_k(\x)^T(\x-P_k(\x)) &\geq f_k(\x) - f_k(P_k(\x)) \nonumber\\
&\geq \frac{\mu}{2}\norm{\x-P_k(x)}^2. \label{gradqg}
\end{align}
Likewise, the following gradient bounds follow from (\textbf{A3})-(\textbf{A4}):
\begin{align}
\norm{\nabla f_k(\x_k)} &= \norm{\nabla f_k(\x_k) - \nabla f_k(P_k(\x_k))} \nonumber\\
&\leq L\norm{\x_k-P_k(\x_k)} = L~\dist{\x_k,\Xk}\label{gradient_norm} \\
\norm{\nabla f_k(\x_k)}^2 &\leq 2L(f_k(\x_k)-f_k^\star)).\label{gradient_square}
\end{align}
\colb{Note that $f_k^\star=f_k(P_k(\x_k))$ and this term will be used interchangeably throughout the analysis. }The following inequality will be utilized towards obtaining the bounds in Sec. V. For positive scalars $u$, $v$, and $w$ that satisfy $u^2 > v$, it holds that
\begin{align}
\sqrt{u^2-v+w^2} &\leq \sqrt{u^2 - v + v^2/4u^2 +  w^2} \nonumber\\
&= \sqrt{u^2\left(1-\frac{v}{2u^2}\right)^2 + w^2} \nonumber\\
&\hspace{0cm}\leq u\left(1-\frac{v}{2u^2}\right) + w\nonumber
\\
& = u - \frac{v}{2u} + w.\label{main_inequality}
\end{align}  
\section{Proof of Lemma \ref{lema1}, \colb{Corollary \ref{p1iidcor}, and Lemma 3}}\label{proof_lem1}
\begin{IEEEproof}[\textbf{Proof of \eqref{contraction}}]
The proof begins by expanding the expression for $\norm{\x_{k+1}-P_k(\x_k)}^2$ in terms of $\dist{\x_k,\Xk}$. Specifically, since ${\mathcal{P}_{\mathcal{X}}}(P_k(\x_k)) = P_k(\x_k)$, it follows from the non-expansiveness property of the projection ${\mathcal{P}_{\mathcal{X}}}(\cdot)$ that, 
\begin{align}
&\norm{\x_{k+1}-P_k(\x_k)}^2	\leq \norm{\x_{k} -\alpha (\nabla f_k(\x_{k}) + \mathbf{e}_{k})- P_k(\x_{k})}^2\nonumber\\
	&=\norm{\x_{k} - P_k(\x_{k})}^2 - 2\alpha \nabla f_k(\x_{k})^T (\x_{k} - P_k(\x_{k}))
	\nonumber
	\\
	&\hspace{10mm} +\alpha^2 \norm{\nabla f_k(\x_{k}) }^2 + \alpha^2\norm{\e_{k}}^2\nonumber
	\\
	&\hspace{10mm}-2\alpha \e^T_{k} (\x_{k}- P_k(\x_{k}))
	+ 2\alpha^2 \e^T_k\nabla f_k(\x_{k}) \label{proof_lemm1_2}
	\\
	&\leq\norm{\x_{k} - {P_k(\x_{k})}}^2 - 2\alpha(f_k(\x_k)-f_k^\star)
	\nonumber
	\\
	&\hspace{10mm} + \alpha^2\norm{\e_{k}}^2 +\alpha^2 \norm{\nabla f_k(\x_{k}) }^2 \nonumber
	\\
	&\hspace{10mm}
	-2\alpha \e_k^T (\x_{k}- {P_k(\x_{k})}) + 2\alpha^2 \e_k^T\nabla f_k(\x_{k}) \label{proof_lemm1_20}
	\end{align}
where the inequality in \eqref{proof_lemm1_20} follows from \eqref{gradqg}. 

The next step is different, depending on the gradient error model in effect. In the general case (cf. \textbf{(A1)}), the last two terms in \eqref{proof_lemm1_20} can be bounded from the use of the Cauchy-Schwarz inequality as
\vspace{0mm}
\begin{align}
&-2\alpha \e_k^T (\x_{k}- {P_k(\x_{k})}) + 2\alpha^2 \e_k^T\nabla f_k(\x_{k}) \nonumber\\
&\qquad\leq 2\alpha\norm{\e_k}\left(\norm{\x_k-P_k(\x_k)} + \alpha\norm{\nabla f_k(\x_k)}\right).
\end{align}
Taking the conditional expectation $\Exc{\cdot}:=\Ex{\cdot \mid \Fk}$, we obtain
\begin{align}
&\Exc{-2\alpha \e_k^T (\x_k- {P_k(\x_{k})}) + 2\alpha^2 \e_k^T\nabla f_k(\x_k)} \label{p1genkey}\\
&\leq 2\alpha(\varepsilon_k + \nu\norm{\nabla f_k(\x_k)})\left(\norm{\x_k-P_k(\x_k)} + \alpha\norm{\nabla f_k(\x_k)}\right). \nonumber
\end{align}
On the other hand, when the errors are zero-mean and i.i.d. (cf. (\textbf{A2})), we have that $\x_k$ is independent of $\e_k$, implying that 
\begin{align}
\Exc{-2\alpha \e_k^T (\x_{k}- {P_k(\x_k)}) + 2\alpha^2 \e_k^T\nabla f_k(\x_k)} &= 0 \label{p1iidkey}
\end{align}
The inequalities in \eqref{p1genkey}-\eqref{p1iidkey} can be encoded into a single inequality by making use of an indicator variable $1_d$ that takes the value 1 when (\textbf{A1}) is in effect and zero when both (\textbf{A1}) and (\textbf{A2}) are in effect. Taking conditional expectation in \eqref{proof_lemm1_20}, we have from \eqref{p1genkey}-\eqref{p1iidkey} that 
\begin{align}
	&	\Exc{\norm{\x_{k+1} - {P_k(\x_{k})}}^2}\nonumber
	\\
	&\leq \norm{\x_{k} - {P_k(\x_{k})}}^2 - 2\alpha(f_k(\x_k)-f_k^\star) \nonumber
	\\
	&\hspace{3mm}+2\alpha 1_d(\varepsilon_k + \nu\norm{\nabla f_k(\x_k)})\left(\norm{\x_k-P_k(\x_k)} + \alpha\norm{\nabla f_k(\x_k)}\right)	\nonumber
	\\
	&\hspace{3mm} + \alpha^2 \!\norm{\nabla \!f_k(\x_k)}^2\!\! + \! \alpha^2 \Exc{\norm{{\e}_{k}}^2}\label{lemm:proof:appen2:first} \\
	&\leq \norm{\x_{k} - {P_k(\x_{k})}}^2 - 2\alpha(f_k(\x_k)-f_k^\star) \nonumber
	\\
	&\hspace{5mm}+2\alpha 1_d\Big(\varepsilon_k\norm{\x_k-P_k(\x_k)} + \nu\norm{\nabla f_k(\x_k)}\norm{\x_k-P_k(\x_k)} \nonumber\\
& \hspace{10mm}	+ \alpha\varepsilon_k\norm{\nabla f_k(\x_k)} + \nu\alpha\norm{\nabla f_k(\x_k)}^2\Big)	\nonumber
	\\
	&\hspace{15mm} + \alpha^2(1 + \nu^2)\norm{\nabla \!f_k(\x_k)}^2 + \alpha^2 \varepsilon_k^2 
	\end{align}
where we have used the error bound in (\textbf{A1}). Next using the inequality $2uv \leq u^2+v^2$, and the inequalities in \eqref{gradient_norm} and \eqref{gradient_square}, we obtain
\begin{align}
	&	\Exc{\norm{\x_{k+1} - {P_k(\x_{k})}}^2} \leq \norm{\x_{k} - {P_k(\x_{k})}}^2 \nonumber\\
	&\ \ \ \ - 2\alpha(1-L\alpha(1+\nu^2)-2\alpha\nu L1_d)(f_k(\x_k)-f_k^\star) \nonumber
	\\
	&\ \ \ \ +2\alpha 1_d\Big(\varepsilon_k(1\!\!+\!\!\alpha L)\norm{\x_k\!\!-\!\!P_k(\x_k)} + \nu L \norm{\x_k-P_k(\x_k)}^2	\Big)	\nonumber\\
	&\ \ \ \  + \alpha^2 \varepsilon_k^2 \nonumber\\
	& =(1+2\alpha\nu L 1_d) \norm{\x_{k} - {P_k(\x_{k})}}^2 +  \alpha^2 \varepsilon_k^2\nonumber\\
	&\ \ \ \  - 2\alpha(1-L\alpha(1+\nu^2+2\nu 1_d))(f_k(\x_k)-f_k^\star) \nonumber	\\
	&\ \ \ \  +2\alpha 1_d\varepsilon_k(1+\alpha L)\norm{\x_k-P_k(\x_k)} 	 \label{preqg}
\end{align}
where for the second term on the right of \eqref{preqg} to be positive, it is required that $\alpha L (1+\nu^2+2\nu1_d) < 1$. Next, the QG property of $f_k$ implies that
\begin{align}
&	\Exc{\norm{\x_{k+1} - {P_k(\x_{k})}}^2} \leq \ell^2\norm{\x_{k} - {P_k(\x_{k})}}^2 \nonumber\\
&\ \ \ \ \ \ \ \ \ \ \ \ \ +2\alpha 1_d\varepsilon_k(1+\alpha L)\norm{\x_k-P_k(\x_k)} + \alpha^2\varepsilon_k^2
\end{align}
where 
\colb{\begin{align}\label{ell_definition}
	\ell^2:=1-\mu\alpha(1\!-\!L\alpha(1\!+\!\nu^2\!+\!2\nu 1_d))\!+\!2\alpha\nu L 1_d < 1.
	\end{align}} Note that it is always possible to find an appropriate $\alpha$ if $\mu > 2\nu L$. The use of Jensen's inequality yields
\begin{align}
&	\Exc{\norm{\x_{k+1} - {P_k(\x_{k})}}} \leq \sqrt{\Exc{\norm{\x_{k+1} - P_k(\x_k)}}^2} \nonumber\\
& \leq \Big(\ell^2\norm{\x_{k} - {P_k(\x_{k})}}^2 + \alpha^2\varepsilon_k^2 \nonumber\\
&\qquad +2\alpha 1_d\varepsilon_k(1+\alpha L)\norm{\x_k-P_k(\x_k)} \Big)^{1/2}. \label{unif1}
\end{align}
Observe that for the case when $1_d = 0$, the right hand side is upper bounded by $\ell\norm{\x_{k} - {P_k(\x_{k})}} + \alpha \varepsilon_k$ from the triangle inequality. In order to obtain a compact expression for the general case, define  
\colb{\begin{align}\label{zeta-updaet}
\zeta &:= \begin{cases} \alpha & 1_d = 0 \\
\frac{\alpha (1+\alpha L)}{\ell} & 1_d = 1
\end{cases} \nonumber\\
&= \alpha + 1_d\alpha(1+\alpha L-\ell)/\ell
\end{align}}
and observe that {$\zeta > \alpha$} since $\ell < 1$. Therefore, it follows from \eqref{unif1} that
\begin{align}
&\Exc{\norm{\x_{k+1} - P_k(\x_k)}} \leq \Big(\ell^2\norm{\x_{k} - {P_k(\x_{k})}}^2 + {\zeta}^2\varepsilon_k^2 \nonumber\\
&\hspace{2cm}+2\alpha 1_d\varepsilon_k(1+\alpha L)\norm{\x_k-P_k(\x_k)} \Big)^{1/2} \\
&\leq \ell \norm{\x_k - P_k(\x_k)} + {\zeta} \varepsilon_k \label{contract2}
\end{align}
which is the required expression. 
\end{IEEEproof}

\begin{IEEEproof}[\textbf{Proof of \eqref{sumdist}}]
This proof will make use of the following bounds, that hold for any $\x^\star \in \cX^\star$,
\begin{subequations}\label{lbub}
\begin{align}
\norm{\E\x-\E\x^\star} \geq \sigma_{\min}\norm{\x-P(\x)} \label{lower_bound}\\ 
\norm{\E\x-\E\y} \leq \sigma_{\max}\norm{\x-\y} \label{upper_bound}
\end{align}
\end{subequations}
where $\sigma_{\max}:=\sigma_{\max}(\E)$ and $\sigma_{\min}:=\sigma_{\min}(\E)$. 

Using the lower bound in \eqref{lower_bound} and from the triangle inequality, it holds for any sequence $\{\x_k^\star\}_{k=1}^K$ that
\begin{align}
&\sum_{k=1}^K\norm{\x_k - P_k(\x_k)} \leq  \frac{1}{\sigma_{\min}} \sum_{k=1}^K\norm{\E\x_k-\E\x_k^\star} \\
&\leq \frac{1}{\sigma_{\min}}\norm{\E\x_1-\E\x_1^\star} \nonumber\\
&\ \ \ + \frac{1}{\sigma_{\min}}\sum_{k=1}^{K-1} \!\!\left(\norm{\E\x_{k+1}-\E\x_k^\star} + \norm{\E\x_{k+1}^\star\!-\!\E\x_k^\star}\right).
 \end{align}
Next, since all points in $\cX_k^\star$ map to a unique point $\u_k^\star$, it holds that $\E\x_k^\star = \E P_k(\x_k)$. Therefore, the upper bound in \eqref{upper_bound} implies that 
\begin{align}
&\sum_{k=1}^K\norm{\x_k - P_k(\x_k)} \leq \frac{\sigma_{\max}}{\sigma_{\min}}\norm{\x_1-\x_1^\star} \nonumber\\
&\ \ \ \ \ \ + \frac{\sigma_{\max}}{\sigma_{\min}}\sum_{k=1}^{K-1}\!\! \left(\norm{\x_{k+1}-P_k(\x_k)}\! +\! \norm{\x_{k+1}^\star\!-\!\x_k^\star}\right). \label{two}
\end{align}
While the result in \eqref{two} holds for any arbitrary $\{\x_k^\star\}$, of particular interest is a sequence with sublinear path length, as defined in \eqref{pl}. Multiplying both sides of \eqref{two} by $\chi=\frac{\sigma_{\min}}{\sigma_{\max}}$, taking expectations, and substituting the result from \eqref{contract2}, we obtain
\begin{align}
\chi\sum_{k=1}^K&\Ex{\norm{\x_k - P_k(\x_k)}} \leq \norm{\x_1-\x_1^\star} + W_K \nonumber\\
&\ \ \ \ \ \ \ \ \ \ \ \ + \ell \sum_{k=1}^K\Ex{\norm{\x_k-P_k(\x_k)}} + {\zeta} \sum_{k=1}^K \varepsilon_k.
\end{align}
Since $\ell < \chi$, it therefore holds that
\begin{align}
\sum_{k=1}^K\Ex{\norm{\x_k - P_k(\x_k)}} \leq \frac{\norm{\x_1-\x_1^\star} + W_K + {\zeta} E_K}{\chi-\ell}.
\end{align}
\end{IEEEproof}

\begin{IEEEproof}[\textbf{\colb{Proof of Lemma \ref{lem:tracking} and equation \eqref{gradient_bound}}}]
Using the bounds in \eqref{lbub} and the triangle inequality, it follows that
\begin{align}
&\norm{\x_{k+1} - P_{k+1}(\x_{k+1})} \leq \frac{1}{\sigma_{\min}} \norm{\E\x_{k+1} - \E \x_{k+1}^\star} \\
&\qquad\leq \frac{1}{\sigma_{\min}} \left(\norm{\E\x_{k+1}-\E\x_k^\star} + \norm{\E\x_{k+1}^\star - \E\x_k^\star}\right) \\
&\qquad= \frac{1}{\sigma_{\min}}\!\! \left(\norm{\E\x_{k+1}\!-\!\E P_k(\x_k)} \!\!+\!\! \norm{\E\x_{k+1}^\star - \E\x_k^\star}\right) \label{gbeq}\\
&\qquad\leq \frac{1}{\chi}\left(\norm{\x_{k+1}-P_k(\x_k)} + \sigma\right)
\end{align}
where the equality in \eqref{gbeq} holds since $\E\x_k^\star = \E P_k(\x_k)$ and the last inequality makes use of the bounded variation property in (\textbf{A5}). Taking expectation, substituting the result in \eqref{contract2}, {and using the bounds is \eqref{boundedekwk},} we obtain
\begin{align}
\Ex{\norm{\x_{k+1} - P_{k+1}(\x_{k+1})}} \leq {\ell}/{\chi}\Ex{\norm{\x_{k}-P_k(\x_k)}} + \frac{{\zeta} \varepsilon + \sigma}{\chi} \label{tracking_bound}
\end{align}
Recursive use of \eqref{tracking_bound} yields 
\begin{align}
\Ex{\dist{\x_{k+1},\cX^\star_{k+1}}} \leq& \left({\ell}/{\chi}\right)^k\dist{\x_1,\cX_1^\star}\nonumber\\ &+{\left[\frac{(1-({\ell}/{\chi})^k)}{(1-{\ell}/{\chi})}\right]\frac{{\zeta} \varepsilon + \sigma}{\chi}}\label{result_tracking}.
\end{align}
\colb{which is the required result in the statement of Lemma \ref{lem:tracking}. The result in \eqref{gradient_bound} can be obtained by using the upper bound $\frac{\ell}{\chi}<1$ in \eqref{result_tracking}, resulting in
	\begin{align}
	\Ex{\dist{\x_{k+1},\cX^\star_{k+1}}} \leq& \dist{\x_1,\cX_1^\star}+{\frac{{\zeta} \varepsilon + \sigma}{\chi}}
	\end{align} 
which, along with \eqref{gradient_norm}, yields the required result.}
\end{IEEEproof}

\section{Proof of Lemma.~\ref{lam_gen_deter} \colb{and Corollary \ref{cor_gen_stoc}}}\label{proof_deter_general}
As earlier, the proofs of Lemma  \ref{lam_gen_deter} and Corollary \ref{cor_gen_stoc} will be developed in a unified manner. To this end, we will again make use of the indicator function $1_d$ that takes values 0 or 1, depending on the assumption under effect. In particular, the following two cases are considered: (a) the general case under (\textbf{A1}) with $\nu = 0$ for which $1_d = 1$ (\colb{Lemma 2}); and (b) the white noise error case under (\textbf{A1})-(\textbf{A2}) for which $1_d = 0$ \colb{Corollary 2}. 

\begin{IEEEproof}[Proof of Lemma \ref{lam_gen_deter}]
Given a sequence $\{\x_k^\star\}$ satisfying \eqref{pl}, and using the first order convexity property of $f_k$, we obtain
\begin{align}
\norm{\x_{k+1} - \x_k^\star}^2 \nonumber
	\\
	&\hspace{-2.1cm}=\norm{\x_{k} - \x_{k}^{\star}}^2 - 2\alpha \nabla f_k(\x_{k}) ^T (\x_{k} - \x_k^\star)+\alpha^2 \norm{\nabla f_k(\x_k)}^2	\nonumber 	\\
	&\hspace{-2.1cm} \ \ \ -2\alpha \e_k^T (\x_k - \x_k^\star)+ \alpha^2\norm{\e_k}^2
	+ 2\alpha^2 \e_k^T\nabla f_k(\x_k) ^T \\
	&\hspace{-2cm}\leq \norm{\x_{k} - \x_{k}^{\star}}^2 - 2\alpha (f_k(\x_{k})- f_k( \x_{k}^{\star})) + \alpha^2 \norm{\mathbf{e}_{k}}^2\nonumber
	\\&\hspace{-2cm}\ \ + 2\alpha^2 L (f_k(\x_{k})- f_k( \x_{k}^{\star}))   +2\alpha^2 \mathbf{e}_{k}^T \nabla f_k(\x_{k})\nonumber
	\\
	&\hspace{-1.5cm}-2\alpha \mathbf{e}_{k}^T (\x_{k}- \x_{k}^{\star})\label{condi}
	\end{align}
Next, we take conditional expectation given $\Fk$ in \eqref{condi}, and consider the two cases separately. In the general case under (\textbf{A1}) with $\nu = 0$, the last two terms can be bounded as in \eqref{p1genkey} to yield
\begin{align}
&\Exc{2\alpha^2 \e_k^T \nabla f_k(\x_k)-2\alpha \e_k^T (\x_k- \x_k^\star)} \nonumber\\
&\hspace{14mm}\leq 2\alpha\Exc{\norm{\e_k}}\left(\alpha\norm{\nabla f_k(\x_k)} + \norm{\x_k-\x_k^\star}\right) \label{condi2}\\
&\hspace{14mm}\leq 2\alpha\varepsilon_k\left(\alpha\norm{\nabla f_k(\x_k)} + \norm{\x_k-\x_k^\star}\right)\label{a11} \\
&\hspace{14mm}\leq \alpha^2(\varepsilon^2_k + \norm{\nabla f_k(\x_k)}^2) + 2\alpha\varepsilon_k R\label{a12}
\end{align}
where recall that $R:=\text{diam}(\cX)$. {Inequality in \eqref{a12} follows from $2uv\leq u^2+v^2$}. Finally, it follows from \eqref{gradient_square} that
\begin{align}
&\Exc{2\alpha^2 \e_k^T \nabla f_k(\x_k)-2\alpha \e_k^T (\x_k- \x_k^\star)} \nonumber\\
&\hspace{12mm}\leq \alpha^2\varepsilon_k^2 + 2\alpha^2L(f_k(\x_k)-f_k(\x_k^\star)) + 2\alpha\varepsilon_kR.
\end{align} 
In the second case, when gradient errors are zero mean and i.i.d., it can be seen from \eqref{condi2} that the right-hand side is zero. Combining the two cases within \eqref{condi}, it follows that
\begin{align}
&\Exc{\norm{\x_{k+1}-\x_k^\star}^2} \leq \norm{\x_k-\x_k^\star}^2 \nonumber\\
& \ \ \ - 2\alpha (f_k(\x_{k})- f_k( \x_{k}^{\star}))  + \alpha^2 \Exc{\norm{\e_k}^2} + \alpha^2 \norm{\nabla f_k(\x_k)}^2 \nonumber\\
	& \ \ \ + 1_d\left( \alpha^2\varepsilon_k^2 + 2\alpha^2L(f_k(\x_k)-f_k(\x_k^\star)) + 2\alpha\varepsilon_kR\right) \\
&	\leq \norm{\x_k-\x_k^\star}^2  + \alpha^2\varepsilon_k^2 (1_d+1) +2\alpha\varepsilon_kR1_d\nonumber\\
& \ \ \ -2\alpha(1\!-\!\alpha L(1\!+\!\nu^2 \!+\!1_d(1\!-\! \nu^2)))(f_k(\x_k)\!-\!f_k(\x_k^\star)) \label{condi3}
	\end{align}
where we have again used \eqref{gradient_square} and the bounds in (\textbf{A1}) and  for the two cases. It is remarked that the last term in \eqref{condi3} is negative for $\alpha L < 1/(1+\nu^2 +1_d(1- \nu^2))$. Henceforth, denote \colb{\begin{align}
	\xi:=& 2\alpha(1-\alpha L(1+\nu^2 +1_d(1- \nu^2)))\label{lemm2_ini1}\\
	s^2_k :=& \alpha^2\varepsilon_k^2 (1_d+1) +2\alpha\varepsilon_kR 1_d.\label{lemm2_ini2}
	\end{align}} Using  \eqref{main_inequality}, it follows that
\begin{align}
&\!\!\!\!\left[\!\Exc{\!\norm{\x_{k+1} \!-\!\x_k^\star}^2}\right]^{1/2}\!\!\! \leq \norm{\x_k\!-\!\x_k^\star}\! - \!\frac{\xi(f_k(\x_k)\!-\!f_k(\x_k^\star))}{\norm{\x_k\!-\!\x_k^\star}} + s_k \nonumber\\
& \leq \norm{\x_k-\x_k^\star} - \frac{\xi}{R}(f_k(\x_k)-f_k(\x_k^\star)) + s_k.
\end{align}
Taking full expectation and using the Jensen's inequality, it follows that
\begin{align}
&\Ex{\norm{\x_{k+1}-\x_k^\star}} \leq \Ex{\left[\Exc{\norm{\x_{k+1}-\x_k^\star}^2}\right]^{1/2}} \nonumber\\
&\leq \Ex{\norm{\x_k-\x_k^\star}} - \frac{\xi}{R}\Ex{f_k(\x_k)-f_k(\x_k^\star)} + s_k
\end{align}
which is the required result. 
\end{IEEEproof}

\footnotesize
\bibliographystyle{IEEEtran} 
\bibliography{IEEEabrv,ref}

\end{document}